\documentclass[12pt]{article}
\usepackage{amssymb,amsmath,latexsym,amsthm,caption2,yfonts,graphicx}
\usepackage[mathscr]{eucal}

\DeclareMathAlphabet{\mathpzc}{OT1}{pzc}{m}{it}

\begin{document}

\newcommand{\be}{\begin{equation}}
\newcommand{\ee}{\end{equation}}
\newtheorem{lemma}{Lemma}
\newtheorem{theorem}{Theorem}
\newtheorem{corollary}{Corollary}
\newtheorem{definition}{Definition}
\newtheorem{example}{Example}
\newtheorem{proposition}{Proposition}
\newtheorem{condition}{Condition}
\newtheorem{assumption}{Assumption}
\newtheorem{conjecture}{Conjecture}
\newtheorem{problem}{Problem}
\newtheorem{remark}{Remark}
\newcommand{\note}{\noindent {\bf Note:\ }}
\newcommand{\la}{\lambda}
\newcommand{\eps}{\varepsilon}
\newcommand{\vph}{\varphi}
\newcommand{\al}{\alpha}
\newcommand{\bet}{\beta}
\newcommand{\gam}{\gamma}
\newcommand{\kap}{\kappa}
\newcommand{\s}{\sigma}
\newcommand{\sig}{\sigma}
\newcommand{\del}{\delta}
\newcommand{\om}{\omega}

\newcommand{\Gam}{\mathnormal{\Gamma}}
\newcommand{\Del}{\mathnormal{\Delta}}
\newcommand{\Th}{\mathnormal{\Theta}}
\newcommand{\La}{\mathnormal{\Lambda}}
\newcommand{\PI}{\mathnormal{\Pi}}
\newcommand{\Sig}{\mathnormal{\Sigma}}
\newcommand{\Ups}{\mathnormal{\Upsilon}}
\newcommand{\Ph}{\mathnormal{\Phi}}
\newcommand{\Ps}{\mathnormal{\Psi}}
\newcommand{\XI}{\mathnormal{\Xi}}
\newcommand{\Om}{\mathnormal{\Omega}}

\newcommand{\C}{{\mathbb C}}
\newcommand{\D}{{\mathbb D}}
\newcommand{\EE}{{\mathbb E}}
\newcommand{\M}{{\mathbb M}}
\newcommand{\N}{{\mathbb N}}
\newcommand{\PP}{{\mathbb P}}
\newcommand{\Q}{{\mathbb Q}}
\newcommand{\R}{{\mathbb R}}
\newcommand{\U}{{\mathbb U}}
\newcommand{\XX}{{\mathbb X}}
\newcommand{\Z}{\mathbb Z}

\newcommand{\calA}{{\cal A}}
\newcommand{\calB}{{\cal B}}
\newcommand{\calC}{{\cal C}}
\newcommand{\calD}{{\cal D}}
\newcommand{\calE}{{\cal E}}
\newcommand{\calF}{{\cal F}}
\newcommand{\calG}{{\cal G}}
\newcommand{\calH}{{\cal H}}
\newcommand{\calI}{{\cal I}}
\newcommand{\calJ}{{\cal J}}
\newcommand{\calK}{{\cal K}}
\newcommand{\calL}{{\cal L}}
\newcommand{\calM}{{\cal M}}
\newcommand{\calP}{{\cal P}}
\newcommand{\calR}{{\cal R}}
\newcommand{\calS}{{\cal S}}
\newcommand{\calT}{{\cal T}}
\newcommand{\calU}{{\cal U}}
\newcommand{\calV}{{\cal V}}
\newcommand{\calW}{{\cal W}}
\newcommand{\calX}{{\cal X}}
\newcommand{\calY}{{\cal Y}}
\newcommand{\calZ}{{\cal Z}}

\newcommand{\scrA}{\mathscr{A}}
\newcommand{\scrL}{\mathscr{L}}
\newcommand{\LL}{\mathscr{L}}
\newcommand{\scrM}{\mathscr{M}}
\newcommand{\scrS}{\mathscr{S}}
\newcommand{\scrX}{\mathscr{X}}
\newcommand{\X}{\mathscr{X}}
\newcommand{\scrY}{\mathscr{Y}}
\newcommand{\Y}{\mathscr{Y}}
\newcommand{\scrR}{\mathscr{R}}
\newcommand{\scrT}{\mathscr{T}}

\newcommand{\uu}{\text{\boldmath$u$}}
\newcommand{\UU}{\text{\boldmath$U$}}

\newcommand{\frA}{\mathfrak{A}}
\newcommand{\frS}{\mathfrak{S}}

\newcommand{\lan}{\langle}
\newcommand{\ran}{\rangle}
\newcommand{\oo}{\overline}
\newcommand{\skp}{\vspace{\baselineskip}}
\newcommand{\noi}{\noindent}
\newcommand{\supp}{{\mbox{supp}\,}}
\newcommand{\diag}{{\rm diag}}
\newcommand{\trace}{{\rm trace}}
\newcommand{\w}{\wedge}
\newcommand{\lt}{\left}
\newcommand{\rt}{\right}
\newcommand{\pl}{\partial}
\newcommand{\prt}{\partial}
\newcommand{\abs}[1]{\lvert#1\rvert}
\newcommand{\norm}[1]{\lVert#1\rVert}
\newcommand{\mean}[1]{\langle#1\rangle}
\newcommand{\til}{\widetilde}
\newcommand{\wt}{\widetilde}
\newcommand{\wh}{\widehat}
\newcommand{\dist}{{\rm dist}}
\newcommand{\grad}{\nabla}
\newcommand{\cd}{\hspace{-.1em}\cdot\hspace{-.1em}}
\newcommand{\arc}[2]{\text{arc}({#1},{#2})}
\newcommand\bn{{\bf n}}
\newcommand\bfe{{\bf e}}
\newcommand\bm{{\bf m}}
\newcommand\bp{{\bf p}}
\newcommand\br{{\bf r}}
\newcommand\bs{{\bf s}}
\newcommand\bt{{\bf t}}
\newcommand\bN{{\bf N}}
\newcommand\ol{\overline}

\def\thelemma{\arabic{section}.\arabic{lemma}}
\def\thetheorem{\arabic{section}.\arabic{theorem}}
\def\thecorollary{\arabic{section}.\arabic{corollary}}
\def\thedefinition{\arabic{section}.\arabic{definition}}
\def\theexample{\arabic{section}.\arabic{example}}
\def\theproposition{\arabic{section}.\arabic{proposition}}
\def\thecondition{\arabic{section}.\arabic{condition}}
\def\theassumption{\arabic{section}.\arabic{assumption}}
\def\theconjecture{\arabic{section}.\arabic{conjecture}}
\def\theproblem{\arabic{section}.\arabic{problem}}
\def\theremark{\arabic{section}.\arabic{remark}}

\newcommand{\manualnames}[1]{
\def\thelemma{#1.\arabic{lemma}}
\def\thetheorem{#1.\arabic{theorem}}
\def\thecorollary{#1.\arabic{corollary}}
\def\thedefinition{#1.\arabic{definition}}
\def\theexample{#1.\arabic{example}}
\def\theproposition{#1.\arabic{proposition}}
\def\theassumption{#1.\arabic{assumption}}
\def\theremark{#1.\arabic{remark}}
}

\newcommand{\beginsec}{
\setcounter{lemma}{0}
\setcounter{theorem}{0}
\setcounter{corollary}{0}
\setcounter{definition}{0}
\setcounter{example}{0}
\setcounter{proposition}{0}
\setcounter{condition}{0}
\setcounter{assumption}{0}
\setcounter{conjecture}{0}
\setcounter{problem}{0}
\setcounter{remark}{0}
}

\renewcommand{\captionfont}{\small\itshape}

\baselineskip=18pt

\title{\bf MIRROR COUPLINGS AND NEUMANN EIGENFUNCTIONS}

\bigskip
\author{\\
{\bf Rami Atar\thanks{Research partially supported by the Israel
Science Foundation (126/02). }}
\smallskip \\
Department of Electrical Engineering \\ Technion - Israel Institute of Technology\\
Haifa 32000, \ Israel\\
atar@ee.technion.ac.il \\
\medskip \\
{\bf Krzysztof Burdzy\thanks{Research partially
supported by NSF Grant DMS-0600206.}}  \smallskip\\
Department of Mathematics, \
University of Washington \\
Seattle,  WA 98195-4350, \ USA \\
burdzy@math.washington.edu}

\date{July 6, 2006}
\maketitle

\begin{abstract}
We analyze a pair of reflected Brownian motions in a planar domain
$D$, for which the increments of both processes form mirror images
of each other when the processes are not on the boundary. We show
that for $D$ in a class of smooth convex planar domains, the two
processes remain ordered forever, according to a certain partial
order.
 This is used to prove that the second eigenvalue is simple for the
Laplacian with Neumann boundary conditions for the same class of
domains.
\end{abstract}

AMS 2000 Subject classification. Primary: 35J05, Secondary: 60H30.

\bigskip
Keywords: Neumann eigenfunctions, reflected Brownian motion,
couplings

\newpage

\section{Introduction}\label{intro}\beginsec

We will prove that the second eigenvalue for the Laplacian with
Neumann boundary conditions is simple for a class of planar convex
domains. We will also present some geometric properties of the
corresponding eigenfunctions. The main tool that we use is a
coupling of a pair of reflected Brownian motions in the domain, for
which the increments of both processes form mirror images of each
other when both processes are not on the boundary. This coupling,
referred to as a {\it mirror} coupling, has been used before to
study properties of Neumann Laplacian eigenfunctions (see
\cite{banbur}, \cite{burrev} and references therein) and, in
particular, has been used in \cite{ab2} to determine whether the
second eigenvalue is simple. That paper was concerned with ``lip
domains'' defined as follows. A lip  domain is a bounded planar
domain that lies between graphs of two Lipschitz functions with the
Lipschitz constant 1. In particular, it has sharp ``left'' and
``right'' endpoints. The current work complements, in a sense, the
results derived in \cite{ab2}, and shows that the technique based on
mirror couplings is also applicable to a class of smooth planar
domains. The earlier paper \cite{banbur}, that also used couplings
in a similar context, showed that the second Neumann eigenvalue is
simple in a convex planar domain if the domain is sufficiently long,
namely, if the ratio of the diameter to width of the domain is
greater than 3.06. If in addition we assume that the domain has a
line of symmetry, the same conclusion can be reached if the ratio of
the diameter to width of the domain is greater than 1.53 (see
Proposition 2.4 of \cite{banbur}). In the current paper we replace
assumptions on the length to width ratio by a set of conditions
that, in particular, allow us to obtain new results for domains that
are not too long.

The motivation for this article comes from the ``hot spots''
conjecture of J.~Rauch which states that the second Neumann
eigenfunction attains its maximum on the boundary of the domain. The
conjecture does not hold in full generality, see
\cite{burwer,fiber,burduke}. It does hold under a variety of extra
assumptions, see \cite{burrev} for a review of literature. This is
related to the question of eigenvalue simplicity because it is often
easier to analyze a single eigenfunction than a class of
eigenfunctions.
 One technical approach to handle both the hot spots conjecture and
the question of eigenfunction
 simplicity is first to change the problem to the mixed
Neumann-Dirichlet problem by identifying the nodal line for the
second eigenfunction (i.e., the line where the eigenfunction
vanishes). This is easily done in symmetric domains (see
\cite{banbur,jernad,pas}). Thus symmetry greatly simplifies the
analysis of eigenfunctions, and removing symmetry from the
assumptions is one of the main technical goals of this paper. The
present paper is the first part of a project which aims at using
this strategy for proving the hot spots conjecture for domains
that are not necessarily symmetric.

\begin{figure}
\centering
 \includegraphics[width=8cm]{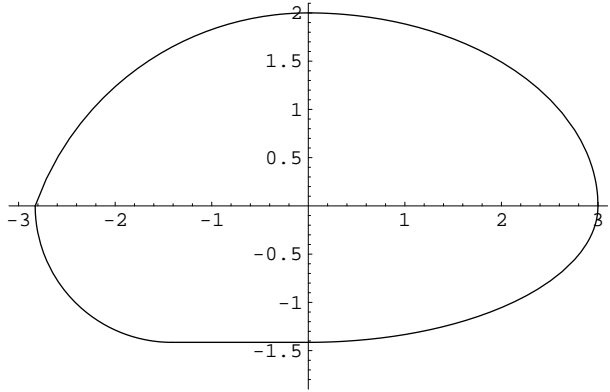}
\caption{A domain with simple second Neumann
eigenvalue.}\label{fige}
\end{figure}

The class of domains that we consider in this paper is defined via a
number of geometric conditions. The conditions are elementary but
their whole set is quite complicated so we will illustrate our main
theorem with some examples. A domain that combines elements of
``extreme'' shapes compatible with our assumptions is depicted in
Fig.~\ref{fige}; see Example \ref{ex1} for the analysis of this
domain.

The set of conditions imposed on a domain $D$ is chosen so that
for appropriately related reflected Brownian motions and an
appropriate partial order, the two processes remain ordered in the
same way forever. We call the line of symmetry for the two
processes a ``mirror.'' We consider mirror couplings, i.e., pairs
of reflected Brownian motions such that the increments of the two
processes are symmetric images of each other with respect to the
mirror, when both processes are in the interior of the domain. The
mirror can be shown to perform a motion that is locally of bounded
variation. The mirror does not move on any interval on which both
processes remain in the interior of the domain. We analyze the
motion of the mirror and construct an appropriate ``Lyapunov
set,'' i.e., a set with the property that the mirror remains in
this set for all times, with probability one, provided that it
starts inside the set.
 The partial ordering alluded to above is defined in terms of this
set. An easy consequence of this property of the coupling is that
there exists a second Neumann eigenfunction that is monotone with
respect to the partial order. We do not know how to prove that the
second eigenvalue is simple using standard results on positive
linear operators such as the Krein-Rutman theorem---to do that, we
would have to impose some extra assumptions on the domain $D$. We
take an alternative approach, similar to that of \cite{ab2}. Along
with the partial order property alluded to above, this approach also
uses crucially the following property of the coupling, which has a
quite complex proof, see \cite{ab2}. If the two processes are
conditioned not to meet up to time 1, the conditional probability
that their distance is greater than $c_1>0$ at time 1 is greater
than $p_1>0$, where $c_1$ and $p_1$ do not depend on the starting
points of the processes.

The paper is organized as follows. In the next section we list
assumptions on the domains that we consider and state our main
result. In Section \ref{secmirror}, we review basic facts about
reflected Brownian motions and mirror couplings. The same section
contains the construction of the Lyapunov set and the proof that
it is left invariant under the dynamics of the mirror process.
Section \ref{secsimple} is devoted to the proof of our main
result, Theorem \ref{th1}. Section \ref{secex} presents some
examples.

We are grateful to Rodrigo Ba\~nuelos for very helpful advice.

\section{Assumptions and the main result}\label{secassn}\beginsec

In the first part of the paper, we consider a bounded strictly
convex planar domain $D$ with $C^2$ boundary $\pl D$. We will
later show that, in a suitable sense, one can remove the
assumptions of strict convexity and $C^2$ smoothness (see the end
of Section \ref{secsimple}). For $A\in\pl D$ let $\bn(A)$ denote
the unit inward normal to $\pl D$ at $A$. For two distinct points
$A,B$ in the plane, we denote by $[A,B]$ the closed line segment
joining them, and by $\ell(A,B)$ the straight line containing
them. We denote by $\scrR[A,B)$ the closed ray contained in
$\ell(A,B)$, starting from $A$ and {\it not} containing $B$. We
fix an orthonormal coordinate system with a basis
$(\bfe_1,\bfe_2)$. We identify $\R^2$ and $\C$ and we use both
types of notation for convenience. For any distinct points $A$ and
$B$, $\angle(A,B)$ denotes the angle between $\bfe_1$ and
$\ell(A,B)$, with the convention that $\angle(A,B) \in [0,\pi)$.
We let $\bp(A,B)=e^{i\angle(A,B)}$, and $\bm(A,B) = -i\bp(A,B)$.
If $\ell$ is a line, we define $\angle\ell$, $\bp(\ell)$ and
$\bm(\ell)$ by choosing any distinct points $A,B\in \ell$ and
letting $\angle\ell= \angle(A,B)$, $\bp(\ell) = \bp(A,B)$ and
$\bm(\ell) = \bm(A,B)$. Note that $\bm(A,B)\cd \bfe_1=\bp(A,B)\cd
\bfe_2\ge 0$. For a point $A\in \prt D$, we let $\angle(A) \in
[0,2\pi)$ be defined by $\bn(A)=e^{i\angle(A)}$.

The closed arc of $\pl D$ joining points $A$ and $B$ on the boundary
is denoted by $\arc{A}{B}$. When we use this notation, we will
specify which one of the two arcs is meant unless it is clear from
the context.

We now list our assumptions on the domain $D$. The assumptions
that are most significant are labelled for future reference in the
proofs.

We will use four sequences of points on the boundary: $P_1,P_2,
\dots, P_6$, $Q_1, Q_2, \dots, Q_6$, $P'_1,P'_2, \dots, P'_6$, and
$Q'_1, Q'_2, \dots, Q'_6$. In this section, we will only discuss
points with subscripts $1,3,4$ and $6$. This is because we chose the
notation so that each of these sequences is naturally ordered along
the boundary, but the existence of points with subscripts $2$ and
$5$ and some special properties will be proved only in Section
\ref{secmirror}.

\begin{figure}
\centering $ \begin{array}{cc}
 \includegraphics[width=8cm]{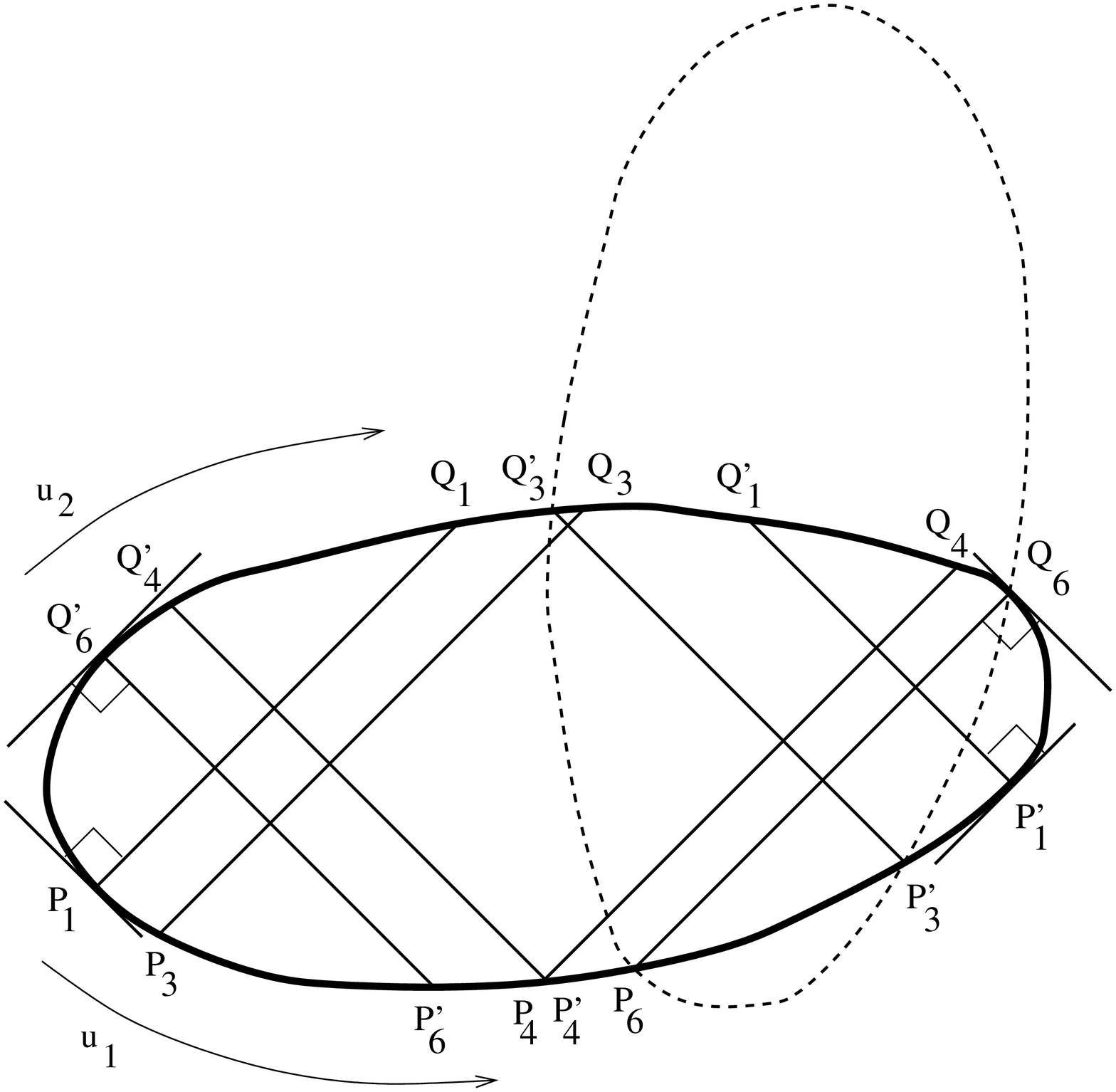}\quad & \quad
 \includegraphics[width=7cm]{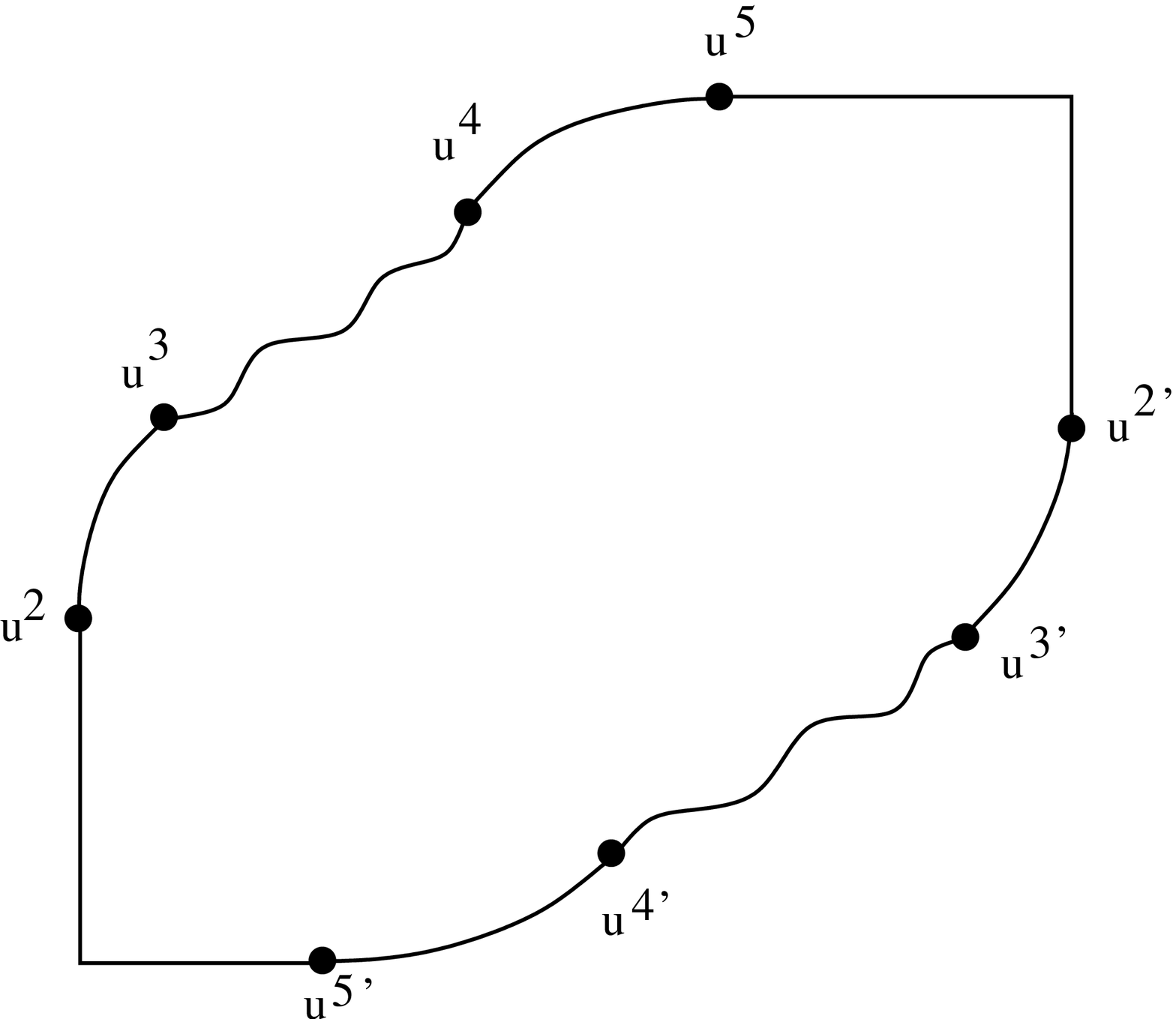}\\
 (a) & (b)
\end{array}
$ \caption{(a) The domain $D$ is shown along with the special points
on its boundary. A reflection of $D$ about $(P'_3,Q'_3)$ is shown in
dashed line. (b) A sketch of the set $\LL$.}\label{figa}
\end{figure}

We assume that there exists an angle $\al\in(0,\pi/2)$ such that all
of the following conditions hold. Let $P_1\in \prt D$ be such that
$\bn(P_1) = e^{i\al}$. Note that $P_1$ exists and is unique because
$D$ is assumed to be strictly convex and $C^2$. Let $Q_1\ne P_1$ be
the unique point on the boundary for which $\angle(P_1,Q_1)=\al$
(see Figure \ref{figa}(a)). Similarly, let $Q_6\in \prt D$ denote
the unique point with $\bn(Q_1) = e^{-i\al}$ and $P_6\in\pl D$ be
such that $\angle(P_6,Q_6)=\al$. We assume that
$(P_6-P_1)\cd\bfe_1>0$ and $(Q_6-Q_1)\cd \bfe_1>0$. We let $\al' =
\pi-\al$ and define points $P'_1$, $Q'_1$, $P'_6$ and $Q'_6$
relative to $\al'$ in the same way that $P_1$, $Q_1$, $P_6$ and
$Q_6$ have been defined relative to $\al$, and assume that
$(P'_6-P'_1)\cd \bfe_1<0$.

Denote by $\prt_\uparrow D$ the closed arc of the boundary from
$Q'_6$ to $Q_6$, not containing $P_1$. We refer to this arc as the
{\it upper part of the boundary}. The arc $\arc{P_1}{P'_1}$ not
containing $Q_6$ will be denoted $\prt_\downarrow D$ and referred
to as the {\it lower part of the boundary}. For points $A,B\in
\prt D$ we write $A<B$ if the first coordinate of $A$ is less than
that of $B$. This ordering will only be used when both $A$ and $B$
are in $\prt_\uparrow D$ or when they are both in $\prt_\downarrow
D$.

We say that a line $\ell$, or line segment $[A,B]$, is {\it
admissible} if it intersects both $\prt_\uparrow D$ and
$\prt_\downarrow D$, and $\angle\ell \in [\al, \al']$ (or
$\angle(A,B) \in [\al, \al']$). For a line $\ell$ that is not
horizontal, we say that a point $C\notin \ell$ is {\it on the
left} of $\ell$ if there exist $D\in \ell$  and $a>0$ such that $C
+ a \bfe_1 = D$. We say that a point is on the left of a line
segment $[A,B]$ if it is on the left of $\ell(A,B)$. Points {\it
on the right} are defined in an analogous way. Suppose $\ell$ is a
line passing through $D$. We say that a boundary point $x\in\pl
D\setminus \ell$ is {\it active for $\ell$} if its reflection
about $\ell$ is in $\ol D$. This seemingly strange term refers to
mirror couplings defined in the next section.

We will state a number of assumptions for $P_1,P_2, \dots, P_6$
and $Q_1, Q_2, \dots, Q_6$. When we say that ``an analogous
condition holds for the primes'' we mean that the analogous
condition holds for $P'_1,P'_2, \dots, P'_6$ and $Q'_1, Q'_2,
\dots, Q'_6$.

\begin{assumption}\label{a1}
There exist line segments $[P_3,Q_3]$ and $[P_4,Q_4]$ satisfying
$\angle(P_3,Q_3)= \angle(P_4,Q_4) =\al$ and such that
$P_1<P_3<P_4<P_6$. Moreover, if $[P,Q]$ is an admissible line
segment with $P_1<P<P_3$ and $\angle(P,Q) \ge \angle (P)$ then no
right boundary point is active. If $[P,Q]$ is an admissible line
segment with $Q_4<Q<Q_6$ and $\angle(P,Q) \ge -\angle (Q)$ then no
left boundary point is active. Analogous conditions hold for the
primes.
\end{assumption}

Suppose that $\ell$ is a line that intersects $D$ and $A\in \prt D
\setminus \ell$ is an active point. Let $\scrT$ denote the line
tangential to $\pl D$ at $A$. If an intersection point of $\ell$
and $\scrT$ exists, it is said to be the {\it hinge} of $A$ at
$\ell$ and it is denoted $H(A,\ell)$. If $\ell=\ell(P,Q)$ then
$H(A,\ell)$ will be called the hinge of $A$ at $[P,Q]$. The name
comes from the fact that the mirror $\ell$ for the coupling of
reflected Brownian motions moves around $H(A,\ell)$ if one of
these processes reflects at $A$ (see Section \ref{secmirror}). We
say that ``hinge $H(A,\ell)$ does not exist'' if $A$ is not an
active point or $\ell$ and $\scrT$ are parallel.

If $P\in \prt_\downarrow D$, $Q\in \prt_\uparrow D$ and $H(A,
\ell(P,Q)) \in \scrR[Q,P)$ then we say that the hinge is {\it
upper}. Otherwise we say that it is {\it lower}. We say that $H(A,
\ell(P,Q))$ is an {\it upper right hinge} if $A$ is on the right
of $[P,Q]$ and $H(A, \ell(P,Q))$ is an upper hinge. We define
upper left, lower right and lower left hinges in an analogous way.

\begin{assumption}\label{a2}
There is $\nu>0$ such that for all $P\in \prt_\downarrow D$ and
$Q\in \prt_\uparrow D$ with $\angle(P,Q)\in[\al-\nu,\al]$ and
$P_3<P<P_4$, there exists no lower left and no upper right hinge.
An analogous condition is assumed for the primes.
\end{assumption}

It follows from Assumption \ref{a3} below that $\arc{P_3}{P_4}$
is, in fact, the largest arc with the above property.

Since $D$ is strictly convex, $\al<\angle(P)$ for $P_1<P<P_3$. We
define $\calA(P_1,P_3)$ as the set of line segments $[P,Q]$ with
the properties $P_1< P< P_3$ and $\angle(P,Q)\in(\al,\angle(P))$.
We define analogously $\calA(Q_4,Q_6)$, $\calA(P'_1,P'_3)$ and
$\calA(Q'_4,Q'_6)$.

It is easy to see that for any $[P,Q] \in\calA(P_1,P_3)$ there
exists at least one lower right hinge. In fact, every
$A\in\prt_\downarrow D$, $A> P$, that is sufficiently close to $P$
is active and the corresponding hinge is lower right.

\begin{assumption}\label{a3}
For any line segment in $\calA(P_1,P_3)$ there exists at least one
lower left but no upper right hinge. For any line segment in
$\calA(Q_4,Q_6)$ there exists at least one upper right but no
lower left hinge. Analogous conditions hold for the primes.
\end{assumption}

For an admissible $[P,Q]$, denote the right [resp., left] part of
the boundary, excluding the endpoints $P$ and $Q$, by $\pl_R(P,Q)$
[$\pl_L(P,Q)$], and its reflection about $\ell(P,Q)$ by
$\til\pl_R(P,Q)$ [$\til\pl_L(P,Q)$].

\begin{assumption}\label{a4}
For every $[P,Q]\in\calA(P_1,P_3)$ [resp., $\calA(Q_4,Q_6)$], the
curves $\pl_R(P,Q)$ and $\til\pl_L(P,Q)$ intersect at a unique
point, and the intersection is nontangential. Moreover, both tangent
lines to these curves at the point of intersection intersect
$\scrR[P,Q)$ [resp., $\scrR[Q,P)$]. Analogous conditions hold for
the primes.
\end{assumption}

\begin{figure}
\centering
 \includegraphics[width=6cm]{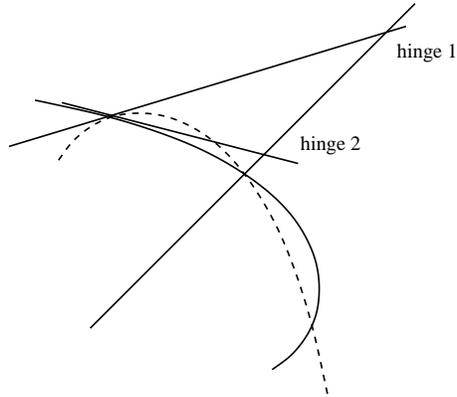}
\caption{A mirror and two hinges.}\label{figb}
\end{figure}

Figure \ref{figb} illustrates a nontangential intersection of the
boundary $\pl D$ and its reflection.

\begin{assumption}\label{a5}
If $[P,Q] \in \calA(P_1,P_3)$ and $[P',Q'] \in \calA(Q'_4Q'_6)$
then $\ell(P,Q) \cap \ell(P',Q')$ is non-empty and belongs to $\ol
D$. An analogous statement holds for the pair $\calA(P'_1,P'_3)$
and $\calA(Q_4,Q_6)$.
\end{assumption}

Our main result is as follows.

\begin{theorem}\label{th1}
Assume that the set $D$ satisfies all the conditions listed in
this section, in particular, Assumptions \ref{a1}-\ref{a5}. Then
the second eigenvalue for the Laplacian in $D$ with Neumann
boundary conditions is simple.
\end{theorem}

\section{Mirror coupling analysis}\label{secmirror}\beginsec

We start by a review of definitions and results from \cite{ab2} on
mirror couplings of reflected Brownian motions.

Let $W$ denote a standard planar Brownian motion and suppose that
$x\in \ol D$. The equation
$$
X(t)=x+W(t)+\int_0^t \bn(X(s))d\hat L(s),
$$
where $\hat L$ denotes the local time of $X$ on the boundary, has
a unique strong solution, referred to as a {\it reflected Brownian
motion}. The local time does not increase when $X$ is away from
the boundary of $D$, i.e., $\int_0^\infty 1_{\{X_s \in D\}} d \hat
L_s = 0$, a.s. By a {\it coupling} we mean a pair of processes
defined on the same probability space. We define a {\it mirror}
coupling of reflected Brownian motions, denoted by $X$ and $Y$ and
starting from $x,y\in\ol D$, by means of the following set of
equations:
\begin{align}\label{10}
X(t)&=x+W(t)+L(t), \quad L(t)=\int_0^t \bn(X(s))d\hat L(s),\\
\label{10+}
Y(t)&=y+Z(t)+M(t),\quad M(t)=\int_0^t \bn(Y(s))d\hat
M(s),
\end{align}
\begin{equation}\label{11}
Z(t)=W(t)-2\int_0^t \bm(s)\,\bm(s)\cd dW(s),\qquad
\bm(t)=\frac{Y(t)-X(t)}{\norm{Y(t)-X(t)}}.
\end{equation}
Here $\hat M$ stands for the local time of $Y$ on $\prt D$. The
definition of $\bm$ given above is different from the meaning
given to this symbol in the previous section but the two vectors
will be effectively identified in our arguments so no confusion
should arise. The equations (\ref{10})-(\ref{11}) have a unique
strong solution up to the time $\zeta=\inf\{t:\lim_{s\to
t-}(X(s)-Y(s))=0\}$ (see \cite{ab2} for the precise meaning of
this statement). The random variable $\zeta$ is called the {\it
coupling time}. While $\{X(t), t\ge0\}$ is well defined by
\eqref{10}, the process $\bm$, and consequently $Y$ is only
well-defined on $[0,\zeta)$. We set
\begin{equation}\label{11+}
\text{ $Y(t)=X(t)$ for all $t\ge\zeta$.} \end{equation}
 Each of the processes $\{X(t), t\ge0\}$ and $\{Y(t),t\ge0\}$ is a
reflected Brownian motion in $D$, and the pair $(X,Y)$ is a strong
Markov process (cf.\ \cite{ab2}).

So long as the processes $X$ and $Y$ have not coupled (i.e., for
times $t<\zeta$), one can talk of a process $\ell(t)$, taking values
in the set of lines in the plane and referred to as the {\it mirror}
process, defined at time $t$ as the line with respect to which
$X(t)$ and $Y(t)$ are symmetric. Clearly, $\bm(t)$ is a unit vector
perpendicular to the mirror, by \eqref{11}. It is also clear that,
for each $t<\zeta$, $X(t)$ and $Y(t)$ can not lie at any boundary
point that is not active for $\ell(t)$. The main result of this
section states that under the assumptions of Section \ref{secassn}
there is a nontrivial subset of $\ol D\times\ol D$ that is left
invariant under the dynamics of the pair $(X,Y)$. It is more
convenient to state and prove this result in terms of the motion of
the mirror $\ell(t)$, a process that is locally of bounded
variation.

We next develop an equation for the intersection points of the
mirror with the boundary. Let $P(t)$ and $Q(t)$ denote the two
intersection points of the mirror $\ell(t)$ with $\pl D$ (for
$t<\zeta$). Let $\bp(t)=\norm{Q(t)-P(t)}^{-1}(Q(t)-P(t))$ and note
that $\bp(t)$ is orthogonal to $\bm(t)$. We label the points in
$\ell(t) \cap \prt D$ in such a way that $\bp(t)=i\bm(t)$. Recall
how $\bp(\ell)$ and $\bm(\ell)$ have been defined in Section
\ref{secassn}. If
 \be\label{51}(Q(t)-P(t))\cd \bfe_2>0,\ee
both definitions of $\bp$ and $\bm$ are consistent in the sense
that $\bp(t) = \bp(\ell(t))$ and $\bm(t) = \bm(\ell(t))$. This is
the case, in particular, when $P(t)\in \prt_\downarrow D$ and
$Q(t) \in \prt_\uparrow D$. Note that by convexity one has
\begin{equation}\label{42}
\bp(t)\cd \bn(P(t))>0,\quad \bp(t)\cd \bn(Q(t))<0.
\end{equation}
It will be convenient to work with the arclength parametrization
of the boundary. If $A\in \prt_\downarrow D$ then we denote by
$U_1(A)$ the length of the arc from $P_1$ to $A$ within
$\prt_\downarrow D$. Analogously, if $A\in \prt_\uparrow D$ then
we write $U_2(A)$ to denote the length of the arc from $Q'_6$ to
$A$ within $\prt_\uparrow D$. We will write $U_1(t) = U_1(P(t))$
and $U_2(Q(t))$ if $P(t)\in \prt_\downarrow D$ and $Q(t)\in
\prt_\uparrow D$. Let $\zeta_0 = \zeta\w\inf\{t\in[0,\zeta):
P(t)\notin \prt_\downarrow D \hbox{ or } Q(t)\notin \prt_\uparrow
D\}$. The process $\{U(t)=(U_1(t),U_2(t)), 0\leq t<\zeta_0\}$
uniquely identifies the mirror process $\ell(t)$ for
$t\in[0,\zeta_0)$. Denote
\begin{equation}\label{100}
V(t)=\norm{X(t)-Y(t)}, \quad \theta(t)=\angle(\ell(t)),\quad
\hbox{ for  } t<\zeta.
\end{equation}
We will suppress the dependence on $t$ for all quantities in the
following lemma.

\begin{lemma}\label{lem5}
We have
\begin{equation}\label{20}
dU_1=(\bp\cd \bn(P)V)^{-1}(-(X-P)\cd dL+(Y-P)\cd dM),
\end{equation}
\begin{equation}\label{21}
dU_2=(\bp\cd \bn(Q)V)^{-1}((X-Q)\cd dL-(Y-Q)\cd dM),
\end{equation}
and \be\label{15} d\theta=V^{-1}\bp\cd(dM-dL), \ee on the time
interval $[0,\zeta_0)$.
\end{lemma}

\noi {\it Remark.} Let
\begin{equation}\label{80}
F=(-(\bp\cd \bn(P)V)^{-1}(X-P)\cd \bn(X), (\bp\cd
\bn(Q)V)^{-1}(X-Q)\cd \bn(X)),
\end{equation}
for $t\in[0,\zeta_0)$ for which $X(t)\in\pl D$ (in which case
$\bn(X)$ is well defined), and set $F=0$ otherwise. Similarly, let
\begin{equation}\label{81}
G=((\bp\cd \bn(P)V)^{-1}(Y-P)\cd \bn(Y),-(\bp\cd
\bn(Q)V)^{-1}(Y-Q)\cd \bn(Y)),
\end{equation}
for $t$ such that $Y\in\pl D$ and $G=0$ otherwise. We can write
equations \eqref{20}- \eqref{21} in the form
\begin{equation}\label{79}
dU=Fd|L|+Gd|M|.
\end{equation}

\proof By the results of \cite{ab2}, the process $\bm$ satisfies
$$
d\bm=V^{-1}(dM-dL)-V^{-1}\bm\,[\bm\cd(dM-dL)],
$$
that can be written as
 \be\label{14} d\bm=V^{-1}\bp\,[\bp\cd(dM-dL)]. \ee
Fix any $t_0\geq 0$ and assume that $\{t_0<\zeta_0\}$ holds. Let
$\bn_0=\bn(P(t_0))$ and $\br_0=-i\bn_0$. Let $Z(t)$ be the
intersection of $\ell(t)$ and the line tangential to $\prt D$ at
$P(t_0)$. Set $x_1(t)=\br_0\cd (Z(t)-P(t_0))$ and $m_1(t)=\br_0\cd
\bm(t)$. It follows from \eqref{14} that
 \be\label{kb1} dm_1=V^{-1}[\bp\cd \br_0]\,[\bp\cd(dM-dL)]. \ee
Elementary geometry can be used to check that
$$
x_1=\frac{(X+Y-2P(t_0))\cd \bm}{2m_1}.
$$
Applying Ito's formula to this representation of $x_1$ yields
\be\label{kb2}
dx_1=(2m_1)^{-1}\bm\cd(dM+dL)+(2m_1)^{-1}(X+Y-2P(t_0))\cd
d\bm-\frac12 m_1^{-2} [\bm\cd(X+Y-2P(t_0))]dm_1. \ee
 Consider any vector $\bs$. Since $\bp \cd \br_0 = - \bn_0 \cd
\bm$, we have
$$
- \bs \cd \bm ( \bp \cd \br_0 + \bn_0 \cd \bm) =0.
$$
We obtain in succession,
$$ \bs\cd[-(\bp \cd \br_0)\bm - (\bn_0 \cd \bm)\bm] =0,
$$
$$ \bs\cd[m_1\bp-(\bp \cd \br_0)\bm -
(\bn_0 \cd \bm)\bm - (\bn_0 \cd \bp)\bp] =0,
$$
$$ \bs\cd[m_1\bp-(\bp \cd \br_0)\bm - \bn_0 ] =0,
$$
$$m_1 \bs \cd \bp - [\bm \cd \bs]\,[\bp \cd \br_0] = \bn_0 \cd \bs.
$$
This, \eqref{14} and \eqref{kb1} imply that
$$
m_1 \bs \cd d\bm - \bm \cd \bs dm_1 = V^{-1} [\bn_0 \cd \bs] \,
[\bp\cd(dM-dL)].
$$
Next we substitute $\bs = X+Y - 2 P(t_0)$ to obtain
$$
m_1 [X+Y - 2 P(t_0)] \cd d\bm - \bm \cd [X+Y - 2 P(t_0)] dm_1 =
V^{-1} [\bn_0 \cd (X+Y - 2 P(t_0))] \, [\bp\cd(dM-dL)],
$$
and
\begin{align}\notag
(2m_1)^{-1} [X+Y - 2 P(t_0)] \cd d\bm &- \frac12 m_1^{-2} \bm \cd
[X+Y - 2 P(t_0)] dm_1 \\
\notag &= (2m_1^2V)^{-1} [\bn_0 \cd (X+Y - 2 P(t_0))] \,
[\bp\cd(dM-dL)].
\end{align}
We combine this with \eqref{kb2} to see that
 \be \label{kb4}
 dx_1=(2m_1)^{-1}\bm\cd(dM+dL)+(2m_1^2V)^{-1}[\bn_0\cd(X+Y-2P(t_0))]\,[\bp\cd
 (dM-dL)].
 \ee
Note that $m_1=\bp\cd \bn_0$ and at time $t=t_0$, the vector $\bp$
is a positive multiple of $(Y+X)/2 - P$. Hence, for $t=t_0$,
$$
m_1^{-1} [\bn_0 \cd ((Y+X)/2 - P)]\, \bp - ((Y+X)/2 - P) =0.
$$
We obtain the following sequence of identities for $t=t_0$,
$$
(Y-X)/2 + m_1^{-1} [\bn_0 \cd ((Y+X)/2 - P)]\, \bp - (Y - P) =0,
$$
$$
V\bm/2 + m_1^{-1} [\bn_0 \cd ((Y+X)/2 - P)]\, \bp - (Y - P) =0,
$$
$$
(2m_1)^{-1}\bm + (2m_1^2V)^{-1} [\bn_0 \cd (Y+X - 2P)]\, \bp -
(m_1V)^{-1}(Y - P) =0,
$$
$$
(2m_1)^{-1}\bm \cd dM + (2m_1^2V)^{-1} [\bn_0 \cd (Y+X - 2P)]\,
\bp \cd dM  = (m_1V)^{-1}(Y - P) \cd dM,
$$
 \be \label{kb3}
 (2m_1)^{-1}\bm \cd dM + (2m_1^2V)^{-1} [\bn_0 \cd (Y+X -
 2P)]\, \bp \cd dM  = ([\bp \cd \bn_0]\,V)^{-1}(Y - P) \cd dM.
 \ee
An analogous calculation yields
 \be \label{kb5}
 (2m_1)^{-1}\bm \cd dL - (2m_1^2V)^{-1} [\bn_0 \cd (Y+X -
 2P)]\, \bp \cd dL  = -([\bp \cd \bn_0]\,V)^{-1}(X - P) \cd dL.
 \ee
We combine \eqref{kb4}-\eqref{kb5} to obtain for $t=t_0$,
$$
dx_1=([\bp\cd \bn(P)]\,V)^{-1}[-(X-P)\cd dL+(Y-P)\cd dM].
$$
The processes $x_1$ and $U_1$ satisfy $\frac{dU_1}{dx_1}(t_0)=1$
because the boundary of $D$ is $C^2$. Therefore \eqref{20}
follows. The proof of \eqref{21} is analogous.

Finally, from $\bp(t)=i\bm(t)=e^{i\theta(t)}$ it is easily seen
that $d\theta=\bp\cd d\bm$, hence by \eqref{14} we obtain
\eqref{15}.
 \qed

\noi{\it Construction of the Lyapunov set}

We will construct a subset of the state space for mirrors
(straight lines in the plane) with the property that if it
contains $\ell(t)$ then it contains $\ell(s)$ for all $s\geq t$,
a.s. It is convenient to encode mirror positions using their
intersection points with $\prt D$ and arclength parametrization
$U_1$ and $U_2$, and so we will work with the process
$U=(U_1,U_2)$ and a set $\LL\subset \R^2$ in the state space of
$U$. Going back to the assumptions and terminology of Section
\ref{secassn}, if $\ell$ is an admissible line, let $P$ and $Q$
denote its intersection points with $\prt_\downarrow D$ and
$\prt_\uparrow D$, and let $u_1=U_1(P)$ and $u_2=U_2(Q)$. Let $\ol
u_1$ and $\ol u_2$ denote the length of $\prt_\downarrow D$ and
$\prt_\uparrow D$. Then for $i=1,2$, $u_i$ takes values in $[0,\ol
u_i]$. We will define the set $\LL$ as a subset of $\UU=[0,\ol
u_1] \times[0,\ol u_2]$. A point $u\in\LL$ thus represents an
admissible line $\ell$. The one-to-one (not onto) map from
admissible line segments to points in $\UU$ described above is
denoted by $\phi$, i.e., the image of $[P,Q]$ is $\phi(P,Q)$. We
use the notation $u^k = \phi([P_k, Q_k])$ and $u^k{}' =
\phi([P'_k, Q'_k])$ for $k=1,2,\dots,6$, where the special line
segments with subscripts $1,3,4$ and $6$ were defined in Section
\ref{secassn}, and those with subscripts $2$ and $5$ will be
defined below.

The boundary of $\LL$ consists of several pieces which will be
described one by one. First, the following set will be a part of
the boundary:
$$
\arc{u^3}{u^4} := \{\phi(P,Q):P_3\le P\le P_4, \angle(P,Q)=\al\}.
$$
Note that this subset of $\UU$ is a curve connecting the points
$u^3$ and $u^4$, corresponding to $[P_3,Q_3]$ and $[P_4,Q_4]$. See
Figure \ref{figa}(b).

Next we describe a curve that begins at the point $u^4$. To this
end we will need the following lemma.

For $A\in \ell(P,Q)\setminus D$, let
$$
d_{P,Q}(A):=\begin{cases} \norm{A-Q}, & A\in \scrR[Q,P),\\
\norm{A-P}, & A\in\scrR[P,Q).\end{cases}
$$
\begin{lemma}\label{lem4}
For every $[P,Q]\in \calA(P_1,P_3)$, there exists a left boundary
point $P_\leftarrow=P_\leftarrow(P,Q)$ with a lower hinge
$H_\leftarrow=H_\leftarrow(P,Q)$ such that $d_{P,Q}(H_\leftarrow)\le
d_{P,Q}(H)$ for every lower left hinge $H$, and there exists a right
boundary point $P_\rightarrow$ with a lower hinge $H_\rightarrow$
such that $d_{P,Q} (H_\rightarrow) \ge d_{P,Q}(H)$ for every lower
right hinge $H$. We also have
 \be\label{kb31}
 d_{P,Q}(H_\leftarrow) > d_{P,Q}(H_\rightarrow).
 \ee
For $[P,Q]\in \calA(Q_4,Q_6)$ there exist points $Q_\leftarrow$
and $Q_\rightarrow$ with properties analogous to $P_\leftarrow$
and $P_\rightarrow$.

Furthermore, for every $\eps>0$ there exists $C_\eps<\infty$ with
the following properties. Suppose that $[P,Q],[\til P,\til
Q]\in\calA(P_1,P_3)$ are such that $\angle(P,Q)<\angle(P)-\eps$,
and assume that a similar inequality holds for $[\til P,\til Q]$.
Then
 \begin{align}\label{kb6}
& \norm{P_\leftarrow- \til P_\leftarrow}+\norm{P_\rightarrow- \til
P_\rightarrow} \le C_\eps(\norm{P-\til P}+\norm{Q-\til Q}),
\end{align}
where $P_\leftarrow=P_\leftarrow(P,Q)$ and $\til
P_\leftarrow=P_\leftarrow(\til P,\til Q)$, etc. Similarly, if
$[P,Q],[\til P,\til Q]\in\calA(Q_4,Q_6)$,
$\angle(P,Q)<-\angle(Q)-\eps$, and a similar inequality holds for
$[\til P,\til Q]$ then
 \begin{align}\label{kb6b}
 & \norm{Q_\leftarrow- \til
Q_\leftarrow}+\norm{Q_\rightarrow- \til Q_\rightarrow} \le
C_\eps(\norm{P-\til P}+\norm{Q-\til Q}).
\end{align}
 Analogous results hold for the primes.
\end{lemma}

\proof Let $[P,Q]\in\calA(P_1,P_3)$. Assumption \ref{a4} asserts the
existence of a unique point of intersection of $\pl_R(P,Q)$ and
$\til\pl_L(P,Q)$. Let $P_\rightarrow\in\pl_R(P,Q)$ denote this
point, and let $P_\leftarrow(P,Q)\in\pl_L(P,Q)$ denote its
reflection about $[P,Q]$. By Assumption \ref{a4}, $P_\rightarrow$
has a lower right hinge, denoted by $H_\rightarrow$, and
$P_\leftarrow$ has a lower left hinge, $H_\leftarrow$. Assumption
\ref{a4} implies that all active points having lower right hinges
must lie on the $\arc{P}{P_\rightarrow}$. Thus the inequality
$d_{P,Q}(H_\rightarrow)\ge d_{P,Q}(H)$ for lower right hinges $H$
follows from convexity. Moreover, no active point having a lower
left hinge can lie on the $\arc{P_\leftarrow}{P}$ (excluding
$P_\leftarrow$), and thus by convexity, $d_{P,Q}(H_\leftarrow)\le
d_{P,Q}(H)$ for lower left hinges $H$. Since by Assumption \ref{a4}
the intersection is nontangential, inequality \eqref{kb31} follows.
Finally, let $\eps>0$ be given. For all line segments
$[P,Q]\in\calA(P_1,P_3)$ satisfying $\angle(P,Q)<\angle(P)-\eps$,
the intersection of $\pl_R(P,Q)$ and $\til\pl_L(P,Q)$ is
nontangential, with a lower bound on the angle of intersection.
Hence by smoothness of $\pl D$, the dependence of the point of
intersection on $P$ and on $Q$ in this class is Lipschitz, with a
constant depending only on $\eps$. It follows from this and the
definition of $P_\rightarrow$ and $P_\leftarrow$ that these two
points are Lipschitz functions of $P$ and on $Q$, with the Lipschitz
constant depending only on $\eps$.
 \qed

For $u\in\UU$ let $[ P(u), Q(u)]$ denote the corresponding line
segment with $ P(u)\in \prt_\downarrow D$, and with an abuse of
notation, let $\bp(u)=\bp( P(u), Q(u))=e^{i\angle(P(u),Q(u))}$. Let
$Q_\leftarrow(u)$ and $Q_\rightarrow(u)$ denote the boundary points
defined relative to $[P(u), Q(u)]$ in Lemma \ref{lem4} above. Note
that $Q_\leftarrow(u)$ has an upper left hinge and
$Q_\rightarrow(u)$ has an upper right hinge. We will prove in Lemma
\ref{lem1} below existence and some properties of a constant $a^*>0$
and a curve $\{u(a), a\in(0,a^*)\}$ in $\UU$ defined by the initial
condition $u(0+)=u^4$ and the following set of ordinary differential
equations,
\begin{equation}\label{60}
\frac d{d a}u_1=\dot u_1=(\bp(u)\cd \bn( P(u)))^{-1}
[-(Q_\leftarrow(u)- P(u))\cd \bn(Q_\leftarrow(u))
-(Q_\rightarrow(u)- P(u))\cd \bn(Q_\rightarrow(u))],
\end{equation}
\begin{equation}\label{61}
\frac d{d a}u_2= \dot u_2=(\bp(u)\cd \bn( Q(u)))^{-1}
[(Q_\leftarrow(u)- Q(u))\cd \bn(Q_\leftarrow(u))
+(Q_\rightarrow(u)- Q(u))\cd \bn(Q_\rightarrow(u))].
\end{equation}
These equations are obtained from \eqref{20}-\eqref{21} by formally
replacing $U$ by $u$, $d|L|$ by $da$, $d|M|$ by $-da$, $V$ by 1,
$X(t)$ by $Q_\leftarrow(u)$ and $Y(t)$ by $Q_\rightarrow(u)$. We
note that we could have formally replaced $d|L|$ by $c_1 da$ and
$d|M|$ by $-c_2 da$; that would not substantially alter the rest of
the argument. The right hand sides of \eqref{60}-\eqref{61} are well
defined for $u\in \phi(\calA(Q_4,Q_6))$ by Lemma \ref{lem4}. The
number $a^*$ has the property that the line $\phi^{-1}(u(a))$ is
asymptotically normal to $\pl D$ as $a\uparrow a^*$; see below for a
precise statement.
 We denote by $u^5$ the limit $\lim_{a\uparrow a^*}u(a)$ (that
exists by the result below) and denote $P_5$, $Q_5$ accordingly.

\begin{lemma}\label{lem1}
There exists a unique constant $a^*\in(0,\infty)$ with the following
properties. The equations \eqref{60}-\eqref{61} have a unique
solution on $(0,a^*)$, and $u(a)\in \phi(\calA(Q_4,Q_6))$ on this
interval. The limit $u^5=\lim_{a\uparrow a^*}u(a)$ exists and one
has $Q_4<Q_5<Q_6$, where $[P_5,Q_5] = \phi^{-1}(u^5)$. Also,
$\lim_{a\uparrow a^*} \bp(u(a))\cd \bn( Q(u(a)))=-1$, i.e., the line
$\phi^{-1}(u(a))$ is asymptotically normal to $\pl D$ at $Q_5$.
Finally, the right hand sides of \eqref{60}-\eqref{61} are positive
on $(0,a^*)$.
\end{lemma}

\proof By convexity of $D$, it follows that $\bp\cd \bn(P)>0$ and
$\bp\cd \bn(Q)<0$ for $[P,Q]\in\calA(Q_4,Q_6)$. For the same reason,
 \begin{equation}\label{900}
 (Q_\leftarrow(u)- P(u))\cd \bn(Q_\leftarrow(u)) <0,
 \qquad
 (Q_\rightarrow(u)- P(u))\cd \bn(Q_\rightarrow(u)) <0,
 \end{equation}
 $$
 (Q_\leftarrow(u)- Q(u))\cd \bn(Q_\leftarrow(u)) <0,\qquad
 (Q_\rightarrow(u)- Q(u))\cd \bn(Q_\rightarrow(u))<0.
 $$
This shows that the right hand sides of \eqref{60} and \eqref{61}
are strictly positive for $[P,Q]\in\calA(Q_4,Q_6)$. Moreover, using
the definition of $\calA(Q_4,Q_6)$, one can see that the left hand
side of the first inequality in \eqref{900} is bounded away from
zero. As a result, the right hand side of \eqref{60} is bounded away
from zero for $[P,Q]\in\calA(Q_4,Q_6)$.

Let $\til Q\in\prt D$ be the point with $\bn(\til Q)=ie^{i\al}$. By
Assumptions \ref{a2} and \ref{a3} there are small perturbations of
$[P_4,Q_4]$ for which there is no upper right hinge, and there are
some for which there exists an upper right hinge. It is not hard to
see that this implies that $Q_\rightarrow(u)\to \til Q$ as $u\to
u^4$ along every sequence for which the hinge exists. We use this to
extend the definition of $Q_\rightarrow(u)$, so that
$Q_\rightarrow(u^4)=\til Q$. Consequently, the right hand sides of
\eqref{60} and \eqref{61} are extended continuously to
$\phi(\calA(Q_4,Q_6)\cup\{[P_4,Q_4]\})$. Let $\calA_\eps$ denote the
set of line segments in $\calA(Q_4,Q_6)$ having
$\angle(P,Q)<-\angle(Q)-\eps$. The local Lipschitz property asserted
in \eqref{kb6}-\eqref{kb6b} and the smoothness of $\bn(\cdot)$
implies that the right hand sides of \eqref{60} and \eqref{61} are
Lipschitz functions of $u$ for
$u\in\phi(\calA_\eps\cup\{[P_4,Q_4]\})$ (with a constant depending
on $\eps$). Let $a_\eps=\inf\{a>0:u(a)\not\in\phi(\calA_\eps)\}$.
The last assertion implies that for every $\eps>0$ there exists a
unique solution on an interval $[0,a_\eps)$, with the initial
condition $u(0)=u^4$. Since by construction
$\angle(P_4,Q_4)<-\angle(Q_4)$, and because the right hand side are
strictly positive, we have that $a_\eps>0$ for all small $\eps>0$.
Since $u_1$ is bounded for $[P,Q]\in\calA(Q_6,Q_6)$, it follows from
the remark above regarding the right hand side of \eqref{60} being
bounded away from zero, that $a_\eps$ are bounded by a finite
constant. The constants $a_\eps$ are clearly monotone, the limit
$a^*=\lim_{\eps\to0}a_\eps$ exists and is finite. The solution to
\eqref{60}-\eqref{61} on $[0,a^*)$ is thus well-defined and unique.
We have already shown that the right hand sides of \eqref{60} and
\eqref{61} are positive. Hence, $u_1$ and $u_2$ are monotone
functions of $a$ and it follows that the limit $u^5:=\lim_{a\uparrow
a^*}u(a)$ exists. We let $[P_5,Q_5]=\phi^{-1}(u^5)$.

We will show that $Q_5< Q_6$. It follows from \eqref{60}-\eqref{61}
that
\begin{equation}\label{106}
\frac{du_1}{du_2}= -\frac{\bp\cd \bn(Q)}{\bp\cd \bn( P)}
 \cdot \frac {(Q_\leftarrow - P)\cd \bn(Q_\leftarrow)
 +(Q_\rightarrow- P)\cd \bn(Q_\rightarrow)}
 {(Q_\leftarrow- Q)\cd \bn(Q_\leftarrow)
 +(Q_\rightarrow- Q)\cd \bn(Q_\rightarrow)} .
\end{equation}
We can consider this as an equation for $u_1$ as a function of $u_2$
with the initial condition $u_1|_{u_2=u_2^4}=u_1^4$. For comparison,
we consider a curve $v(a)=v=(v_1,v_2)$ in $\UU$ with the property
that $\angle(\phi^{-1}(v))=\al$ for all $a>0$. This curve satisfies
$$
\frac{dv_1}{dv_2}=-\frac{\bp\cd \bn( Q)}{\bp\cd \bn( P)}.
$$
Again, we can consider the above as an equation for $v_1$ as a
function of $v_2$, with the same initial condition as for $u_1$,
namely, $v_1|_{v_2=u_2^4}=u_1^4$ (this is because
$\angle(\phi^{-1}(u^4))=\al$). The fact that $Q_\rightarrow$ has an
upper hinge implies that $(Q-P)\cd \bn(Q_\rightarrow)<0$. Similarly,
$(Q-P)\cd \bn(Q_\leftarrow)<0$. It follows that the second fraction
on the right hand side of \eqref{106} is strictly less than 1.
Standard comparison results for univariate ODE's imply that
$u_2<v_2$ whenever $u_1=v_1$. This shows that
\begin{equation}\label{901}
\text{$\angle(\phi^{-1}(u(a)))>\al=\angle(P_6,Q_6)$ for every
$a\in(0,a^*]$.}
\end{equation}
We are in the middle of an argument that is supposed to show that
$Q_5< Q_6$. We now argue by contradiction and assume that $Q_5\ge
Q_6$. Then $[P,Q_6]=\phi^{-1}(u(\hat a))$ for some $\hat
a\in(0,a^*]$ and $P$. By \eqref{901}, $\angle(\phi^{-1}(u(\hat
a))>\al=-\angle(Q(u(\hat a)))$.
 Hence for small $\eps>0$, $\angle(\phi^{-1}(u(\til
a))>-\angle(Q(u(\til a)))$ for an appropriate $\til a< a_\eps$.
This contradicts the definition of $a_\eps$. We conclude that
$Q_5<Q_6$.

Finally, note that the limit $\lim_{a\uparrow a^*} \bp(u(a))\cd \bn(
Q(u(a)))$ exists by monotonicity of $u$ and is equal to
$\bp(u(a^*))\cd \bn( Q(u(a^*)))$. Since by \eqref{901} we have
$\angle(P_5,Q_5)>\al$, and since $Q_4<Q_5<Q_6$, it follows from the
definitions of $\calA(Q_4,Q_6)$, $\calA_\eps$ and $a_\eps$ that for
all small $\eps>0$,
$\angle(\phi^{-1}(u(a_\eps)))=-\angle(Q(u(a_\eps)))-\eps$. Thus
$\angle(P_5,Q_5)=-\angle(Q_5)$ i.e., $\bp(u(a^*))\cd \bn(
Q(u(a^*)))=-1$. \qed

The part of the boundary constructed above is denoted by
$\arc{u^4}{u^5}$.

Analogously to $\arc{u^3}{u^4}$, we construct
$\arc{u^3{}'}{u^4{}'}$. Similarly to $\arc{u^4}{u^5}$, we
construct $\arc{u^4{}'}{u^5{}'}$, and then $\arc{u^2}{u^3}$ and
$\arc{u^2{}'}{u^3{}'}$.

Next, consider the two line segments $[P_5,Q_5]$ and
$[P'_2,Q'_2]$. Since $Q_5<Q_6$, Assumption \ref{a5} implies that
these two line segments intersect in $D$. As a result, $P_5<P'_2$
and for similar reasons, $Q'_2 < Q_5$. We add the following pieces
to the boundary of $\LL$:
$$
\{\phi(P,Q_5):P_5\le P\le P'_2\}, \quad \{\phi(P'_2,Q): Q'_2\le Q\le
Q_5\}.
$$
We denote this by $\arc{u^5}{u^2{}'}$. Note that it has the form
$$
\{(u_1,u_2)\in\UU:u_1^5\le u_1\le u_1^2{}',
u_2=u_2^5\}\cup\{(u_1,u_2)\in\UU:u_1=u_1^2{}', u_2^2{}'\le u_2\le
u_2^5 \}.
$$
Finally, we construct $\arc{u^2}{u^5{}'}$, the last part of the
boundary of $\LL$, in a way analogous to the construction of
$\arc{u^5}{u^2{}'}$.

In view of Lemma \ref{lem1} it is easy to see that the pieces of
$\prt \LL$ constructed above do not intersect each other, except
for the endpoints. The Lyapunov set $\LL$ is defined as the simply
connected, bounded, closed domain with the boundary comprised of
all arcs constructed above.

\noi{\it Invariance of the set $\LL$.}

Recall the definitions of the mirror coupling $(X(t),Y(t))$, and
$U(t)$, $\ell(t)$ and $\zeta$ from the beginning of this section.
The main result of this section states that the process $U$
remains in $\LL$ if it start in $\LL$.

\begin{theorem}\label{th2} Assume that $X(0)\ne Y(0)$,
$U(0)\in\LL$ and that $X(0)$ is on the left of $\ell(0)$. Then,
with probability 1, for all $t<\zeta$, $U(t)\in\LL$ and
$\bfe_1\cd(Y(t)-X(t))>0$.
\end{theorem}

\proof Suppose for the moment that $U(t)\in\LL$ for all $t<\zeta$.
Then the assertion that $\bfe_1\cd(Y(t)-X(t))>0$ for $t<\zeta$
follows from sample path continuity of $X$ and $Y$ and the fact
that $\ell(t)\in[\al,\pi-\al]$, $t<\zeta$. Hence, it remains to
show that $U(t)\in\LL$ for all $t<\zeta$.

We define $\tau=\inf\{t\in[0,\zeta):U(t)\not\in\LL\}$ and
$E=\{\tau<\zeta\}$, with the convention $\inf\emptyset=+\infty$.
We will show that $\PP(E)=0$. On $E$, we let
$u^*=U(\tau)\in\pl\LL$. Note that $U(t_k)\in\LL^c$ along a
sequence $t_k\to\tau+$.

Consider first the case when $u^*\in \arc{u^4}{u^5}\setminus
\{u^4,u^5\}$. By \eqref{60}-\eqref{61}, the vector
$\bt^*=(\bt^*_1,\bt^*_2)$ given by the following formula is
tangent to $\pl\LL$ at $u^*$,
\begin{equation}\label{102}
\bt^*_1=(\bp^*\cd \bn( P^*))^{-1} [-(Q_\leftarrow^*- P^*)\cd
\bn(Q_\leftarrow^*) -(Q_\rightarrow^*- P^*)\cd
\bn(Q_\rightarrow^*)],
\end{equation}
\begin{equation}\label{103}
\bt^*_2=(\bp^*\cd \bn( Q^*))^{-1} [(Q_\leftarrow^*- Q^*)\cd
\bn(Q_\leftarrow^*) +(Q_\rightarrow^*- Q^*)\cd
\bn(Q_\rightarrow^*)].
\end{equation}
The superscripts ${}^*$ in the above formula indicate that all the
functions are evaluated at $u^*$. Note that $V^*>0$ since
otherwise we would have $\tau=\zeta$. By Lemma \ref{lem1},
$\bt_i^*>0$, $i=1,2$. Thus $\bN^*:=(\bt^*_2,-\bt^*_1)$ is an
inward normal vector to $\pl\LL$ at $u^*$. It is not necessarily
true that $\norm{\bN^*}=1$. Note that $X^*\in\pl D$ or $Y^*\in\pl
D$ (or both) because $X$ and $Y$ are continuous and the mirror
$\ell(t)$ is not moving when the reflected Brownian motions are
inside $D$. In the case when $X^*\in\pl D$, the expression on the
right hand side of \eqref{80} evaluated at $\tau$ will be denoted
by $F^*$. Similarly, in the case $Y^*\in\pl D$, $G^*:=G(\tau)$
(cf.\ \eqref{81}).

We will now show that $F^*\cd \bN^*>0$ in the case $X^*\in\pl D$.
Let
$$ \gamma^*=-([\bp^*\cd \bn(P^*)]\, [\bp^*\cd \bn(Q^*)]\, V^*)^{-1},
$$
and note that $\gamma^* >0$. Since $\bN^*=(\bt^*_2,-\bt^*_1)$,
\begin{align}\label{83}\notag
F^*_1\bt^*_2-F^*_2\bt^*_1= \gamma^*\{ & (X^*-P^*)\cd \bn(X^*)
[(Q_\leftarrow^*- Q^*)\cd \bn(Q_\leftarrow^*) +(Q_\rightarrow^*-
Q^*)\cd \bn(Q_\rightarrow^*)] \\
& -(X^*-Q^*)\cd \bn(X^*)[(Q_\leftarrow^*- P^*)\cd
\bn(Q_\leftarrow^*) +(Q_\rightarrow^*- P^*)\cd
\bn(Q_\rightarrow^*)]\}.
\end{align}

By Assumption \ref{a3} the hinges corresponding to $X^*$ and to
$Q_\leftarrow^*$ are upper; thus $\bp^*\cd \bn(X^*)<0$ and
$\bp^*\cd \bn(Q_\leftarrow^*)<0$. Also, Lemma \ref{lem4} states
that the distance from the hinge corresponding to $Q_\leftarrow^*$
to $Q^*$ is not smaller than the distance from the hinge
corresponding to $X^*$ to $Q^*$. It follows that the distance from
the hinge corresponding to $Q_\leftarrow^*$ to $P^*$ is not
smaller than the distance from the hinge corresponding to $X^*$ to
$P^*$. One can express this fact by the following inequality:
$$
\frac{(Q_\leftarrow^*-P^*)\cd \bn(Q_\leftarrow^*)}{\bp^*\cd
\bn(Q_\leftarrow^*)}\ge \frac{(X^*-P^*)\cd \bn(X^*)}{\bp^*\cd
\bn(X^*)}.
$$
Since $Q^*-P^*$ is a positive multiple of $\bp^*$, it follows that
$$
[(Q_\leftarrow^*-P^*)\cd \bn(Q_\leftarrow^*)]\,[(Q^*-P^*)\cd
\bn(X^*)] \geq [(X^*-P^*)\cd \bn(X^*)]\, [(Q^*-P^*)\cd
\bn(Q_\leftarrow^*)]
$$
and
\begin{align*}
&[(Q_\leftarrow^*-P^*)\cd \bn(Q_\leftarrow^*)]\,[(Q^*-P^*)\cd
\bn(X^*) - (X^*-P^*)\cd \bn(X^*)]\\
&\ge  [(X^*-P^*)\cd \bn(X^*)]\, [(Q^*-P^*)\cd \bn(Q_\leftarrow^*)
- (Q_\leftarrow^*-P^*)\cd \bn(Q_\leftarrow^*)].
\end{align*}
This gives
\begin{equation}\label{85}
[(X^*-P^*)\cd \bn(X^*)]\,[(Q_\leftarrow^*-Q^*)\cd
\bn(Q_\leftarrow^*)] -[(X^*-Q^*)\cd
\bn(X^*)]\,[(Q_\leftarrow^*-P^*)\cd \bn(Q_\leftarrow^*)]\ge0.
\end{equation}
Next, by Assumption \ref{a3} and Lemma \ref{lem4} (applied to
$A(Q_4Q_6)$), $Q_\rightarrow^*$ has an upper right hinge. By Lemma
\ref{lem4}, the distance of this hinge from $Q^*$ is strictly
larger than that of the hinge corresponding to $Q_\leftarrow^*$,
and, in turn, that corresponding to $X^*$. It follows that the
distance of the hinge corresponding to $Q_\rightarrow^*$ from
$P^*$ is strictly larger than that of the hinge corresponding to
$X^*$. This can be written as
\begin{equation}\label{104}
\frac{(Q_\rightarrow^*-P^*)\cd \bn(Q_\rightarrow^*)}{\bp^*\cd
\bn(Q_\rightarrow^*)} > \frac{(X^*-P^*)\cd \bn(X^*)}{\bp^*\cd
\bn(X^*)}.
\end{equation}
A calculation similar to the one leading to \eqref{85} yields the
strict inequality
\begin{equation}\label{86}
[(X^*-P^*)\cd \bn(X^*)]\,[(Q_\rightarrow^*-Q^*)\cd
\bn(Q_\rightarrow^*)] -[(X^*-Q^*)\cd
\bn(X^*)]\,[(Q_\rightarrow^*-P^*)\cd \bn(Q_\rightarrow^*)]
> 0.
\end{equation}
We add \eqref{85} and \eqref{86} and combine the result with
\eqref{83} to obtain $F^*\cd \bN^*>0$.

A similar calculation (that is slightly more complicated due to
the fact that $Y^*$ can either have an upper or a lower hinge)
results in the conclusion that $G^*\cd \bN^*>0$ in the case
$Y^*\in\pl D$.

We now go back to \eqref{80}-\eqref{79}. By the sample path
continuity of the processes $|L|$, $|M|$, $X$, $Y$, $\bp$, $P$,
$Q$ and the continuity of the vector field $\bn$ on $\pl D$, the
fact that $F^*\cd \bN^*>0$ provided that $X^*\in\pl D$ implies
that there exists a (random) $\eps>0$ such that for all $\del>0$
small enough,
$$
\int_{[\tau,\tau+\del]}[F(t)\cd \bN^*]\, d|L|(t)
> \eps|L|([\tau,\tau+\del]).
$$
Similarly, for sufficiently small $\del>0$,
\begin{equation}\label{kb8}
\int_{[\tau,\tau+\del]}[G(t)\cd \bN^*] \, d|M|(t)
> \eps|M|([\tau,\tau+\del]).
\end{equation}
Thus by \eqref{79}, for all sufficiently small $\del>0$,
\begin{equation}\label{84}
(U(\tau+\del)-U(\tau))\cd \bN^*\ge\eps\la(\del),
\end{equation}
where $\la(s)=|L|([\tau,\tau+s])+|M|([\tau,\tau+s])$.

The boundary of $\LL$ is $C^1$ in a neighborhood of $u^*$ so for
any sequence $\hat u_k\in \LL^c$ with $\hat u_k\to u^*$, we have
$\limsup_{k\to \infty} (\hat u_k-u^*)\cd \bN^* / \norm{\hat
u_k-u^*} \leq 0$. Since $U(t_k)\in\LL^c$ for a sequence
$t_k\to\tau+$, we have $\la(s)>0$ for $s>0$. The last two
observations imply that
$$
(U(t_k)-U(\tau))\cd \bN^*\le C\la(t_k)r(t_k),
$$
for some (random) $C<\infty$ and $r(t)$ such that $r(0+)=0$. This
contradicts \eqref{84}; thus the probability that $U$ exits $\LL$
through $u^*\in \arc{u^4}{u^5}\setminus \{u^4,u^5\}$ is equal to
zero.

Next consider the case when $U$ exits $\LL$ through
$\arc{u^5}{u^2{}'}$, excluding the points $u^5$ and $u^2{}'$. By
Assumption \ref{a1}, $X(\tau)$ is not on the boundary. Thus
\begin{equation}\label{902}
|L|([\tau-\eps,\tau+\eps])=0\quad\text{for an appropriate random
$\eps>0$}.
\end{equation}
Hence only the term $Gd|M|$ is present in \eqref{79}. It is easy to
see from \eqref{81} that $G^*_i<0$ for both $i=1,2$. The inward unit
normals to $\pl\LL$ at $u^*$ can be either $(0,-1)$ or $(-1,0)$,
except there is a single point (corner) where any convex combination
of these vectors points inside $\LL$. Thus $G^*\cd\bN^*>0$ in all
cases and \eqref{kb8} holds. The argument following \eqref{kb8} can
be now repeated to rule out the possibility of exiting through
$\arc{u^5}{u^2{}'}$.

Consider now exit through $u^5$ ($u^2{}'$ can be treated similarly).
By Lemma \ref{lem1}, the line $[P(u),Q(u)]$ is asymptotically normal
to $\pl D$ at $Q_5$. This implies that $Q_\leftarrow(u)-Q(u)$ and
$Q_\rightarrow(u)-Q(u)$ vanish as $u\to u^5$ along the curve
$\arc{u^4}{u^5}$. Thus the right hand side of \eqref{61} vanishes in
this limit, and it follows that $\pl\LL$ is $C^1$ at $u^5$, with the
unit inward normal $(0,-1)$ at this point. As in the preceding
paragraph, Assumption \ref{a1} implies \eqref{902}. The analysis of
this case can now be finished by the same argument as in the case of
an exit thorough $\arc{u^5}{u^2{}'}$.

Consider now the possibility that $U$ exits through
$\arc{u^3}{u^4}$, excluding $u^3$ and $u^4$. Recall $\theta$
defined in \eqref{100} and set $\theta^*=\theta(\tau)$. Then
$\theta^*=\al$. If the trajectory of $U$ exits $\LL$ at time
$\tau$ then the trajectory of $\theta$ exits $[\al,\pi-\al]$ at
the same time. Thus for every $\eps>0$ there exist $s,t$ such that
$\tau\le s<t<\tau+\eps$ and
\begin{equation}\label{101}
\theta(r)<\theta(s)=\al,\quad \text{ for all } r\in(s,t].
\end{equation}
Recall from \eqref{15} that $d\theta=V^{-1}\bp\cd(dM-dL)$. By
Assumption \ref{a2}, there is no upper right or lower left hinge
within $(s,t)$, if $\eps$ is small enough. This means a right
hinge is necessarily lower. Thus, within this time interval,
$\bp\cd dM=[\bp\cd \bn(Y)]\, d|M|$ and $\bp\cd \bn(Y)\ge0$.
Similarly, $\bp\cd \bn(X)\le0$. As a result,
$\theta(t)\ge\theta(s)$, contradicting \eqref{101}. We see that
$U$ cannot exit $\LL$ through $\arc{u^3}{u^4}$.

The discussion of the possible exit through $u^4$ will be split
into two steps---one similar to the treatment of $\arc{u^4}{u^5}$
and the second one similar to that of $\arc{u^3}{u^4}$. One can
show that the interior angle formed by $\prt \LL$ at $u^4$ is less
than or equal to $\pi$ but the calculation will not be provided
here. If the angle is greater than $\pi$ then the first step of
the argument given below would alone suffice to complete the
proof.

Let us thus review the argument provided for $\arc{u^4}{u^5}$.
Note that, by Lemma \ref{lem1}, the formula
\eqref{102}-\eqref{103} for the one-sided tangent line to this arc
is still valid for $u^*=u^4$, and that
$\bn(Q_\rightarrow^*)=i\bp^*$ in this case. Consider a closed set
$\LL'$ whose boundary is the same as that of $\LL$, except that
the arc joining $u^3$ and $u^4$ is replaced by a curve that, in
the vicinity of $u^4$, coincides with a ray starting from $u^4$
and is oriented as $-\bt^*$. Note that $\pl\LL'$ is $C^1$ at
$u^4$, and, as before, $\bN^*=(\bt_2^*,-\bt_1^*)$ is an inward
normal to $\LL'$. Recall that we have assumed that
$U(t_k)\not\in\LL$ for a sequence $t_k\to\tau+$. We begin by
showing that $U(t)\in\LL'$ for all $t\in(\tau,\tau+\eps)$, if
$\eps>0$ is small enough. As in the argument for $\arc{u^4}{u^5}$,
we can achieve that by showing that $F^*\cd \bN^*>0$ and $G^*\cd
\bN^*>0$. The argument leading to \eqref{85} holds. The one that
leads to \eqref{86} is not valid since $\bp^*\cd
\bn(Q_\rightarrow^*)=0$ and \eqref{104} can not be used. To obtain
\eqref{86}, note that $(P^*-Q^*)\cd \bn(Q_\rightarrow^*)=0$, and,
therefore, the left hand side of \eqref{86} can be written as
$$
[(Q_\rightarrow^*-P^*)\cd \bn(Q_\rightarrow^*)]\, [(Q^*-P^*)\cd
\bn(X^*)].
$$
We have $(Q_\rightarrow^*-P^*)\cd \bn(Q_\rightarrow^*)<0$, and,
since $X^*$ has an upper left hinge, $(Q^*-P^*)\cd \bn(X^*)<0$.
Thus \eqref{86} holds. The argument following \eqref{86} can be
repeated and one concludes that $U(t)\in \LL'$ for all
$t\in[\tau,\tau+\eps)$, where $\eps>0$ is sufficiently small.

Next, note that $\bt^*\in\R_+^2$. Hence if $B=B(u^4,\rho)$ denotes a
disc and $\calC =  B\cap \{(u_1,u_2)\in\R^2: u_1<u^4_1,u_2<u^4_2\}$
then for sufficiently small $\rho>0$ we have
$B\cap(\LL'\setminus\LL)\subset \calC$. As a result, if
$[P,Q]=\phi^{-1}(u)$ for any $u\in \calC$ then $P_3<P<P_4$ and
$Q_3<Q<Q_4$. Moreover, given any point $u\in \calC$, one can find a
point $u'\in\arc{u^3}{u^4}$ such that $u'_1>u_1$ and $u'_2=u_2$.
Then $\angle \phi^{-1}(u) < \angle \phi^{-1}(u')$. The angle for
each such $u'$ equals $\al$, by construction of $\arc{u^3}{u^4}$.
Thus, by Assumption \ref{a2}, an upper right hinge does not exist
for $\phi^{-1}(u)$. The argument that we used for $\arc{u^3}{u^4}$
can now be adapted to show that $U$ cannot exit $\LL$ through $u^4$.

The proof is analogous for the other parts of the boundary of
$\LL$.
 \qed

\section{Multiplicity of the second eigenvalue}\label{secsimple}

\beginsec
In this section we prove Theorem \ref{th1}. The overall strategy
of the proof is similar to that in \cite{ab2}. We begin by
reformulating our main tool, Theorem \ref{th2}, in a convenient
way.

First, given $(x,y)\in\ol D\times\ol D$, $x\ne y$, let $m(x,y)$ be
the line of symmetry for $x$ and $y$ and let $\{P_{x,y},Q_{x,y}\} =
m(x,y) \cap \prt D$, with the convention that the second coordinate
of $P_{x,y}$ is less than or equal to that of $Q_{x,y}$. Let
$T^1\subset\ol D\times\ol D$ denote the set of pairs $(x,y)$, $x\ne
y$, for which $\phi(P_{x,y},Q_{x,y})\in\LL$. Let
$$
T=\{(x,y)\in T^1: \bfe_1\cd(y-x)>0\}.
$$
For $(x,y)\in T$, let $\PP_{x,y}$ denote a probability measure
under which $(X(t),Y(t))$ is a mirror coupling starting from
$(x,y)$, and recall that $Y(t)=X(t)$ for all $t\ge\zeta$. Let
$\EE_{x,y}$ denote the corresponding expectation. An alternative
statement of Theorem \ref{th2} is that if $(x,y)\in T$ then
$\PP_{x,y}$-a.s., $(X(t),Y(t))\in T$ for all $t<\zeta$.

Let $D_L$ denote the connected component of $D \setminus
\ell(P_1,Q'_6)$ not containing $Q_6$ in its closure, and similarly
let $D_R$ denote the connected component of $D \setminus
\ell(P'_1,Q_6)$ not containing $Q'_6$ in its closure. Because
$\LL$ is constructed as a subset of $\UU$, it follows from Theorem
\ref{th2} that, for $(x,y)\in T$, $\PP_{x,y}$-a.s., for all
$t<\zeta$, $\ell(t)$ does not intersect $D_L$ or $D_R$. Thus, for
$(x,y)\in T$,
\begin{equation}\label{99}
X(t)\not\in \oo{D_R} \text{ and } Y(t)\not\in\oo{D_L}\quad \text{
for all } t<\zeta,\ \PP_{x,y}\text{-a.s.}
\end{equation}

We will use the following well known probabilistic representation
of solutions to the heat equation. Suppose that $f_0$ is a bounded
function on $\ol D$. Let $f:[0,\infty)\times\ol D\to\R$ denote the
solution to $(1/2) \Delta f = (\prt/\prt t) f$ with initial values
$f(0,x)=f_0(x)$ and Neumann boundary conditions on $\pl D$. Then
 \be\label{kb9}
 f(t,x)=\EE_{x} f_0(X(t)).
 \ee
In particular, if $\mu_2>0$ is the second eigenvalue for the
Laplacian with Neumann boundary conditions and $\psi$ is any
eigenfunction corresponding to $\mu_2$ then the above formula may
be applied to $f(t,x)=e^{-\mu_2 t} \psi(x)$ and we obtain
 \be\label{kb10}
 \psi(x)= e^{\mu_2t} \EE_{x} f_0(X(t)).
 \ee

\begin{lemma}\label{lemfour}
There exist constants $c_1,p_1>0$ such that for every $(x,y)\in
T$,
$$
\PP_{x,y}(\norm{X(1)-Y(1)}\ge c_1\mid \zeta>1)\ge p_1.
$$
\end{lemma}

\proof The assertion is the same as in Lemma 4 of \cite{ab2}. The
proof is very similar to that in \cite{ab2} with minor, obvious
adaptations, and is thus omitted. \qed

Let
$$
S=\{f\in C(\ol D):f(y)-f(x)\ge0 \text{ for all } (x,y)\in T\},
$$
$$
\til S=\{f\in C(\ol D):f(y)-f(x)>0 \text{ for all } (x,y)\in T\}.
$$

\begin{lemma}\label{lemf}
If $\psi$ is a second Neumann eigenfunction and $\psi\in S$ then
$\psi\in\til S$.
\end{lemma}
\proof Consider a second Neumann eigenfunction $\psi$ and assume
that $\psi\in S$. Given $(x,y)\in T$, we shall show that
$\psi(y)>\psi(x)$. Let us begin with $(x,y)\in T^o$ (the interior
of $T$). Let $\eps>0$ be so small that $B(x,\eps)\times
B(y,\eps)\subset T^o$. Since $\psi(x')\le\psi(y')$ for $(x',y')\in
T$ and $\psi$ is a non-constant real analytic function on $D$,
there must exist $(x',y')\in B(x,\eps)\times B(y,\eps)$ where
$\psi(x')<\psi(y')$. Thus there also exist balls $B_1\subset
B(x,\eps)$, $B_2\subset B(y,\eps)$ and $\del>0$ such that
$\psi(x')+\del<\psi(y')$ for all $(x',y')\in B_1\times B_2$.
Consider a coupling of processes $(\til X(t),\til Y(t))$, in which
$\til X(t)$ and $\til Y(t)$ are independent Brownian motions
starting from $x$ and $y$, resp., until $\tau:=\inf\{t>0:(\til
X(t),\til Y(t))\in\pl(B(x,\eps)\times B(y,\eps))\}$, at which time
they switch to a mirror coupling. Clearly $\til X(t)$ and $\til
Y(t)$ are reflected Brownian motions in $\ol D$, starting from $x$
and $y$. By Theorem \ref{th2} and the strong Markov property
applied at $\tau$, the process $(\til X(t),\til Y(t))$ does not
leave the set $T$ for $t<\zeta$, a.s. Thus, using \eqref{kb10}, we
obtain
$$
e^{-\mu_2}(\psi(y)-\psi(x))=\EE_{x,y}(\psi(\til Y(1))-\psi(\til
X(1)))\ge\del\PP_{x,y}(H),
$$
where $H$ denotes the event that $\til X$ and, respectively, $\til
Y$ do not leave $B(x,\eps)$ and $B(y,\eps)$ before time 1, and
$\til X(1)\in B_1$, $\til Y(1)\in B_2$. By well known properties
of the standard Brownian motion, the probability of $H$ is
strictly positive. This shows that $\psi(y)>\psi(x)$ for all
$(x,y)\in T^o$.

To complete the proof, it suffices to show that for every
$(x,y)\in T$, a mirror coupling $(X,Y)$ starting from $(x,y)$,
reaches the interior of $T$ by time 1, and $\zeta>1$, with
positive probability. To this end it suffices to show that the
process $U$, if it starts on $\pl\LL$, enters the interior of
$\LL$ before time 1, and $\zeta>1$, with positive probability. We
analyze different parts of $\pl\LL$ separately. If $U(0) \in
\arc{u^3}{u^4}$, consider $z\in\pl D$ and let $z'$ be the mirror
image of $z$ with respect to $[P,Q]=\phi^{-1}(U(0))$. Choose $z$
so that it has an upper left hinge and is located so close to $Q$
that for some $\eps \in(0,\norm{z-z'}/2)$ we have
$B(z',\eps)\subset D$. Let $D'$ be the connected component of $D
\setminus [P,Q]$ that is on the left of $[P,Q]$. Consider the
following event,
\begin{align}\notag
\{&X(t)\in D' \text{ for } t\in[0,1/2]; X(t)\in B(z,\eps) \text{
for } t\in(1/2,1]; \\
\notag & Y(t)\in B(z',\eps) \text{ for } t\in(1/2,1];
|L|([0,1])>0\}.
\end{align}
It is standard to prove that the above event has a strictly
positive probability. Since $B(z',\eps)\subset D$, we have
$\zeta>1$ if this event occurs. Since $\bp\cd \bn(z)<0$, it easily
follows from \eqref{15} that $\theta(1)>\theta(0)=\al$. Thus the
trajectory $U$ enters $\LL^o$ if this event holds.

A similar argument applies for $U(0)\in\arc{u^4}{u^5}$ with $z$
being a point on the boundary, close enough to $P$, having a lower
right hinge (by Assumption \ref{a3} there is no lower left hinge
for $[P,Q]\in \calA(Q_4,Q_6)$ hence a lower right hinge must
exist). Here one uses equations \eqref{20}-\eqref{21} to show that
$U$ enters $\LL^o$. For $U(0)\in\arc{u^5}{u^2{}'}$, take $z$ to be
any boundary point to the right of $(P,Q)$ and use again
\eqref{20} and \eqref{21}.

Finally, consider the special boundary points $u^4$ and $u^5$. Now
that it has been shown that the interior is reached from anywhere
in $\pl\LL$ save these special points, it suffices to show that
the mirror line $\ell(t)$ simply moves (with positive probability)
if it starts at the corresponding positions. However, the only way
it can happen that the mirror does not move with probability 1 is
when the domain $D$ is symmetric with respect to $\ell(0)$. This
is clearly not the case for either $u^4$ or $u^5$, due to our
assumptions. \qed

The following lemma essentially follows from Lemma 4.1 of
\cite{ata}, except that it has slightly weaker smoothness
assumptions. The proof given here is shorter than that in
\cite{ata}.

\begin{lemma}\label{kb20}
If $\psi$ is a Laplacian eigenfunction with Neumann boundary
conditions corresponding to $\mu_2$ in a convex bounded domain $D$
then $sup_{x\in D}\norm{\nabla\psi(x)} < \infty$.
\end{lemma}
\proof Consider any points $x,y\in \ol D$ and let $(X,Y)$ be a
mirror coupling of reflected Brownian motions in $D$. Recall that
$\zeta$ stands for the coupling time of $X$ and $Y$. By
\eqref{kb10},
 $$
 |\psi(y)-\psi(x)|= e^{\mu_2} |\EE_{x,y}(\psi( Y(1))-\psi( X(1)))|
 \leq e^{\mu_2} \norm{\psi}_\infty \PP_{x,y} (\zeta >1).
 $$
Since $D$ is a convex domain, $\norm{\psi}_\infty < \infty$ (see,
e.g., \cite{banbur}). An application of the It\^o formula and
equations \eqref{10}-\eqref{11} show that
$$
\norm{X-Y}=\norm{x-y}+\oo{W}+\oo V,
$$
where $\oo W=-2\int_0^\cdot m\cd dW$ and $\oo V=\int_0^\cdot
(n(Y)\cd md|M|-n(X)\cd md|L|)$. The process $\oo W$ is a one
dimensional Brownian motion (with the diffusion constant different
from the standard one) and, by convexity of the domain, the
process $\oo V$ is non-increasing. Hence,
 $$\PP_{x,y} (\zeta >1)
 \leq \PP \left(\inf_{0\leq t \leq 1} \ol W_t >-\norm{x-y}\right)
 \leq c_1 \norm{x-y}.
 $$
We see that, for some $c_2<\infty$,
 $$|\psi(y)-\psi(x)| \leq e^{\mu_2} \norm{\psi}_\infty c_1 \norm{x-y}
 \leq c_2 \norm{x-y}.
 $$
The lemma follows easily from this bound.
 \qed

For $\eps>0$, let
$$
T_\eps=\{(x,y)\in T:\norm{x-y}\ge\eps\}.
$$

\begin{lemma}\label{lems}
Let $c_1$ be as in Lemma \ref{lemfour}. For every
$\eps_1\in(0,c_1)$ such that the interior of $T_{\eps_1}$ is
non-empty and every $\del,\kap>0$ there exists $\eps_2>0$ with the
following property. If $\psi$ is a second Neumann eigenfunction
satisfying
\begin{equation}\label{90}
\psi(y)-\psi(x)\ge\del\quad \text{ for all } (x,y)\in T_{\eps_1},
\end{equation}
\begin{equation}\label{91}
\psi(y)-\psi(x)\ge0\quad\text{ for all } (x,y)\in T_{\eps_2},
\end{equation}
and
\begin{equation}\label{92}
\norm{\nabla\psi}\le\kap\quad\text{ on } D,
\end{equation}
then $\psi\in S$.
\end{lemma}
\proof Let $c_1$ and $p_1$ be as in the statement of Lemma
\ref{lemfour}. Fix any $\eps_1\in(0,c_1)$ such that the interior
of $T_{\eps_1}$ is non-empty and consider any $\del,\kap>0$. Let
\begin{equation}\label{98}
p_2=\inf_{(x',y')\in T_{c_1}}\PP_{x',y'}((X(1),Y(1))\in
T_{\eps_1}).
\end{equation}
It follows easily from Lemma 2 of \cite{ab2} that $p_2>0$. Set
$\eps_2=\min(\kap^{-1}\del p_1 p_2,\eps_1)$ and note that
$\eps_2>0$. Let $\psi$ be a second Neumann eigenvalue satisfying
\eqref{90}-\eqref{92}. By \eqref{kb10},
$$
e^{-\mu_2t}(\psi(y)-\psi(x))=\EE_{x,y}[\psi(Y(t))-\psi(X(t))],\quad
t\ge0.
$$
Thus, it suffices to show that, for all $(x,y)\in T$,
\begin{equation}\label{97}
\EE_{x,y}[\psi(Y(2))-\psi(X(2))\mid \zeta>1]\ge0.
\end{equation}
To this end, note that, in view of \eqref{92},
$$
\EE_{x,y}[\psi(Y(2))-\psi(X(2))\mid \zeta>1]\ge
\EE_{x,y}[(\psi(Y(2))-\psi(X(2)))
1_{\{\norm{X(2)-Y(2)}\ge\eps_2\}}\mid \zeta>1]-\kap\eps_2.
$$
Since $(X(2),Y(2))\in T$ a.s., the indicator function on the right
hand side of the last formula can be replaced by
$1_{\{(X(2),Y(2))\in T_{\eps_2}\}}$. By \eqref{91}, the above
inequality remains valid if the indicator function is further
replaced by $1_{\{(X(2),Y(2))\in T_{\eps_1}\}}$. Thus by
\eqref{90} and \eqref{98},
\begin{align*}
\EE_{x,y}[\psi(Y(2))-\psi(X(2))\mid \zeta>1]  &\ge
\EE_{x,y}[(\psi(Y(2))-\psi(X(2))) 1_{\{(X(2),Y(2))\in
T_{\eps_1}\}}\mid \zeta>1]-\kap\eps_2\\
&\ge
\del\,\PP_{x,y}((X(2),Y(2))\in T_{\eps_1}\mid \zeta>1)-\kap\eps_2\\
&\ge \del\eta(x,y)p_2-\kap\eps_2,
\end{align*}
where
$$
\eta(x,y)=\PP_{x,y}((X(1),Y(1))\in T_{c_1}\mid \zeta>1).
$$
By Lemma \ref{lemfour}, $\eta(x,y)\ge p_1$. Thus
$$
\EE_{x,y}[\psi(Y(2))-\psi(X(2))|\zeta>1] \ge \del p_1 p_2
-\kap\eps_2\ge0,
$$
and we have shown that \eqref{97} holds for all $(x,y)\in T$. This
completes the proof of the lemma. \qed

\noi{\it Proof of Theorem \ref{th1}.} We start by showing that,
whether the second Neumann eigenvalue is simple or not, there
exists a corresponding eigenfunction that lies in $S$. The
multiplicity of $\mu_2$ is either one or two, see
\cite{banbur,nad1,nad2}. Consider first the case when the
multiplicity of $\mu_2$ is two and let $\psi$ and $\psi'$ be
orthogonal Neumann eigenfunctions corresponding to $\mu_2$ and
normalized so that $\int_D\psi^2=\int_D(\psi')^2=1$. Since $\psi$
is real analytic in $D$, it is impossible that it vanishes on all
of $D_R$; thus assume without loss of generality that it is
strictly positive in some ball $B\subset D_R$. Let $f_0$ be a
continuous nonnegative, nonzero function on $\ol D$, supported
inside $B$. Let $f:[0,\infty)\times\ol D\to\R$ denote the solution
to the heat equation in $\ol D$ with the initial values
$f(0,x)=f_0(x)$ and Neumann boundary conditions on $\pl D$. The
function $f$ has the following eigenfunction expansion:
 \be\label{kb11}
f(t,x)=C_1+(C_2\psi(x)+C'_2\psi'(x))e^{-\mu_2t}+R(t,x),
 \ee
where $C_1$, $C_2$ and $C'_2$ are suitable constants, and
$\lim_{t\to\infty}e^{\mu_2t}\sup_{x\in\bar D}|R(t,x)|=0$ (see
\cite{banbur}, Proposition 2.1). Note that $C_2=\int_D f_0\psi>0$.
Therefore $\psi_0:=C_2\psi+C'_2\psi'$ is a nonzero eigenfunction
corresponding to $\mu_2$. We have by \eqref{kb11} for $(x,y)\in
T$,
$$
e^{\mu_2t}(f(t,y)-f(t,x))=\psi_0(y)-\psi_0(x)+\eps(t,x,y)
$$
where $\eps(t,x,y)\to0$ as $t\to\infty$. We now use \eqref{kb9} to
write
$$
\psi_0(y)-\psi_0(x)=e^{\mu_2t}\EE_{x,y}[f_0(Y(t))-f_0(X(t))]-\eps(t,x,y).
$$
By \eqref{99} and the properties of $f_0$, it follows that
$\psi_0(y)-\psi_0(x)\ge0$ for all $(x,y)\in T$.

In the case when $\mu_2 $ is simple, we take $\psi'\equiv 0$ and
repeat the argument to conclude that $\psi \in S$. Obviously, if
we assume that $\mu_2$ is simple, there is no logical need to
prove any properties of eigenfunctions to finish the proof of
Theorem \ref{th1}. However, the fact that $\psi\in S$ is an
interesting by-product of the proof.

In what follows, $\psi_0$ denotes an eigenfunction in $S$.

To prove that $\mu_2$ is simple, we use a variation of a proof
from \cite{ab2}. We argue by contradiction and assume that $\mu_2$
is not simple. Let $\psi_1$ denote a second Neumann eigenfunction
orthogonal to $\psi_0$. It follows from Lemma \ref{lemf} that
$-\psi_1$ and $\psi_1$ cannot both lie in $S$; we thus assume
without loss of generality that $\psi_1\not\in S$. Let
$$
\psi_a=(1-a)\psi_0+a\psi_1,\quad 0<a<1,
$$
$$
a^*=\inf\{a\in[0,1]:\psi_a\not\in S\}.
$$
We claim that $a^*<1$. If $a^*=0$, then we are done. Otherwise,
for $a<a^*$ and $(x,y)\in T$, one has $\psi_a(y)-\psi_a(x)\ge0$.
By continuity of the function $a\to \psi_a$,
$\psi_{a^*}(y)-\psi_{a^*}(x)\ge0$ for all $(x,y)\in T$. Thus
\begin{equation}\label{95}
\psi_{a^*}\in S,
\end{equation}
and, therefore, $a^*<1$. This implies that
\begin{equation}\label{94}
\exists \  a_k\downarrow a^*,\ a_k\in(a^*,1),\ \psi_{a_k}\not\in
S.
\end{equation}
For $a=a_k$ as above, let
$$
\eps(a)=\sup\{\norm{x-y}:(x,y)\in T,\ \psi_a(y)-\psi_a(x)<0\}.
$$
By \eqref{95} and Lemma \ref{lemf}, $\psi_{a^*}\in\til S$. When
$k\to\infty$, $\psi_{a_k}\to\psi_{a^*}$ uniformly in $\ol D$. This
easily implies that,
\begin{equation}\label{96}
\eps(a_k)\to0 \text{ as } k\to\infty.
\end{equation}
Fix some $\eps_1>0$ as in Lemma \ref{lems}. Since
$\psi_{a^*}\in\til S$, we have $\psi_{a^*}(y)-\psi_{a^*}(x)>0$ for
all $(x,y)$ in the closed set $T_{\eps_1}$. Thus there are
constants $\del>0$ and $k_0$ such that for all $k>k_0$ and
$(x,y)\in T_{\eps_1}$, $\psi_{a_k}(y)-\psi_{a_k}(x)\ge\del$. By
Lemma \ref{kb20}, $\kap := \sup_D(\norm{\nabla\psi_0}
+\norm{\nabla\psi_1})<\infty$. Let $\eps_2>0$ be defined relative
to $\eps_1, \del $ and $\kap$ as in Lemma \ref{lems}. By
\eqref{96}, we have that $\psi_{a_k}(y)-\psi_{a_k}(x)\ge0$ for all
$(x,y)\in T_{\eps_2}$, provided that $k$ is large enough. Thus by
Lemma \ref{lems}, $\psi_{a_k}\in S$ for all large $k$. This
contradicts \eqref{94}. We conclude that $\mu_2$ is simple. \qed

\begin{corollary}\label{kbgrad}
Let $D_M$ be the part of $D$ between $[P_3,Q'_4]$ and $[P'_3, Q_4]$.
Suppose that the assumptions of Theorem \ref{th1} hold. Consider a
second eigenfunction $\psi$ and define $\beta(x)$ by
$\nabla\psi(x)=e^{i\beta(x)}$. Then
$\beta(x)\in[\al-\pi/2,\pi/2-\al]$ for all $x\in D_M$, or this
assertion holds for $-\psi$.
\end{corollary}

\proof The corollary is an easy consequence of the fact that
$\psi\in S$, established in the first part of the proof of Theorem
\ref{th1}.
 \qed

\begin{proposition}\label{kb21}
Suppose that the second eigenvalue $\mu_2$ for the Laplacian with
Neumann boundary conditions in a Lipschitz domain $D\subset \R^d$
is simple and there exist disjoint subsets $D', D''\subset \ol D$
with non-empty interiors, non-empty open balls $B', B'' \subset
D$, and a coupling of reflected Brownian motions $(X,Y)$ in $D$
such that for any $X(0)=x\in B'$ and $Y(0)=y\in B''$ we have $X(t)
\notin D'$ and $Y(t) \notin D''$ for all $t < \zeta := \inf\{t\geq
0: X(t)=Y(t)\}$, a.s. Let $\psi$ be an eigenfunction corresponding
to $\mu_2$. Then $\psi(z) \geq 0$ for all $z\in D'$ and $\psi(z)
\leq 0$ for all $z\in D''$, or this assertion applies to $-\psi$.
\end{proposition}

\proof Let $p_t(\,\cd\,,\,\cd\,)$ denote the transition density for
the reflected Brownian motion in $D$. Consider any $x,y \in \ol D$.
 Then, by Proposition 2.1 of \cite{banbur}, for some $c_1\in \R$,
$c_2,c_3\in(0,\infty)$, depending on $x$ and $y$,
 \be\label{kb22}
 p_t(x,z) - p_t(y,z) =
 c_1 e^{-\mu_2 t} \psi(z) + R(t,z),
 \ee
and $ |R(t, z)|\leq c_2 e^{-(\mu_2+c_3) t}$, for all $t\geq 1$ and
all $z\in D$. Recall that $\psi$ is a real analytic function that is
not identically constant so it is not constant on balls $B'$ and
$B''$. Hence, we can choose $x\in B'$ and $y\in B''$ such that
$\psi(x)\ne \psi(y)$ and, therefore, $c_1= \psi(x) - \psi(y) \ne 0$.
 Assume without loss of generality that $c_1>0$. Consider any
non-empty open ball $B$ in the interior of $D'$. Then
 \begin{align}\notag
 \int_B (p_t(x,z) - p_t(y,z) ) dz
 &= \PP_{x,y} (X(t) \in B) - \PP_{x,y} (Y(t) \in B)\\ \notag
 &= \PP_{x,y} (X(t) \in B; t< \zeta)- \PP_{x,y} (Y(t) \in B; t<\zeta)
 \\ \notag
 &= - \PP_{x,y} (Y(t) \in B; t<\zeta)
 \leq 0.
 \end{align}
This shows that $p_t(x,\,\cd\,) - p_t(y,\,\cd\,)$ is non-positive
in the interior of $D'$. We combine this with \eqref{kb22} and let
$t\to \infty$ to see that $\psi(z)\leq 0$ for $z\in D'$.
Similarly, $\psi(z)\geq 0$ for $z\in D''$. The inequalities are
reversed if $c_1<0$.
 \qed

We note that the coupling of reflected Brownian motions in
Proposition \ref{kb21} is not assumed to be the mirror coupling.
Among currently known couplings, the mirror coupling seems to be
the only one which can satisfy the assumptions of Proposition
\ref{kb21}. However, some new couplings are proposed from time to
time (see, e.g., \cite{ab1,pas}) so the proposition might be
applied in the future to some other class of couplings.

\begin{corollary}\label{kb23}
Suppose that $D$ is as in Theorem \ref{th1} and recall the
definitions of $D_L$ and $D_R$ from the beginning of this section.
The second eigenfunction $\psi$ for the Laplacian with Neumann
boundary conditions in $D$ is non-negative on $D_L$ and
non-positive on $D_R$, or this assertion applies to $-\psi$.
\end{corollary}

\proof The corollary follows from Proposition \ref{kb21} and
\eqref{99}.
 \qed

The geometric location of the nodal line (i.e., zero set of the
second eigenfunction) was studied in \cite{ab1}. The results of
that paper are logically independent from Proposition \ref{kb21}
in the following sense. The techniques developed in \cite{ab1}
cannot be used to prove Corollary \ref{kb23}. On the other hand,
the location of the nodal line in obtuse triangles is determined
with greater accuracy in \cite{ab1} than it could be done using
Proposition \ref{kb21}.

We will next show how one can remove, in a sense, the assumptions of
strict convexity and $C^2$-smoothness from Theorem \ref{th1}.
Suppose that $D\subset \R^2$ is bounded and convex but not
necessarily strictly convex and $\prt D$ is not necessarily
$C^2$-smooth. Suppose that $\{D_k\}_{k\geq 1}$ is a non-decreasing
sequence of domains satisfying assumptions of Theorem \ref{th1} and
converging to $D$, i.e., $\bigcup_{k\geq 1} D_k = D$. Let $P^k_1$
and $Q'_{6,k}$ be the points defined relative to $D_k$ and analogous
to $P_1$ and $Q'_6$ in $\prt D$. Recall that these points are used
to define the arc parametrization for parts of $\prt D_k$. Assume
that there exist $P_1$ and $Q'_6$ such that $P^k_1 \to P_1$ and
$Q'_{6,k} \to Q'_6$ as $k\to\infty$. We make similar assumptions
about existence of points $Q^k_1$ and $P'_{6,k}$ and their
convergence to $Q_1$ and $P'_6$. Let $D_L,D_R\subset D$ be defined
as at the beginning of this section, in terms of $P_1, Q'_6, P'_1$
and $Q_6$ described above. Let $\LL_k$ be the Lyapunov set
corresponding to the domain $D_k$ and let $\LL= \ol{\limsup_{k\to
\infty} \LL_k}$, i.e., $\LL= \ol{ \bigcap_{n\geq 1} \bigcup_{k\geq
n} \LL_k}$. Let $T$ be defined as at the beginning of this section,
relative to the present definition of $\LL$.

\begin{proposition}\label{kb24}
Assume that the above conditions for $D$ hold and every line of
symmetry for $D$ is horizontal or vertical (hence $D$ may have two,
one or no lines of symmetry). Assume that $D$ is not a rectangle,
and that $D_L, D_R$ and $T$ have non-empty interiors.
 Then the assertion of Theorem \ref{th1} holds, i.e., the second
eigenvalue for the Laplacian with Neumann boundary conditions in
$D$ is simple.
\end{proposition}

\proof First we will prove that there exists a mirror coupling of
reflected Brownian motions $(X,Y)$ in $D$ for which $\LL$ is a
Lyapunov set, i.e., if $(x,y)\in T$ and $(X(0),Y(0)) = (x,y)$ then
$(X(t),Y(t))\in T$ a.s., for all $t<\zeta := \inf\{t\geq0:
X(t)=Y(t)\}$.

Fix any $(x,y)\in T$. It follows from the definition of $\LL$ and
$T$ that there exist $x_k,y_k \in D_k$ such that $x_k\to x$, $y_k\to
y$ and $(x_k,y_k) \in T_k$, where $T_k$ is defined relative to
$D_k$. For any $k$, let $(X_k, Y_k)$ be the mirror coupling of
reflected Brownian motions in $D_k$ with $(X_k(0), Y_k(0)) = (x_k,
y_k)$ as defined in \eqref{10}--\eqref{11+}. We construct all these
processes on a single probability space and use the same process $W$
to define $X_k,Y_k,Z_k,\zeta_k$ and $\bm_k$ for all $k$. By Theorem
2.3 of \cite{burch}, $X_k$'s converge in distribution in the local
uniform topology to a reflected Brownian motion in $D$ and an
analogous statement is true for $Y_k$'s. In particular, there exists
a filtration $(F_t)$ with respect to which the Brownian motion $W$
is a martingale, and there exist processes $X$, $Y$, $Z$, $L$ and
$M$ such that $Z$ is a Brownian motion and an $(F_t)$-martingale,
$X$, $Y$, $L$ and $M$ are $(F_t)$-adapted, and the equations
\eqref{10}-\eqref{10+} hold. Passing to a subsequence if necessary,
$X_k$ converge to $X$ and $Y_k$ to $Y$ locally uniformly, with
probability one. Let $\zeta=\inf\{t:\lim_{s\to t-}(X(s)-Y(s))=0\}$.
By the uniform convergence result, for every $\delta>0$ one has
$\zeta_k\ge\zeta-\del$ for all large $k$. In particular, for all
large $k$, $\bm_k$ is well-defined for $t\le\zeta-\del$, and $Z_k$
satisfies
$$
Z_k(t)=W(t)-2\int_0^t\bm_k(s)\bm_k(s)\cd dW(s)
$$
for $t\le\zeta-\del$. The processes $\bm_k(\cd\w(\zeta-\del))$
converge uniformly to $\bm(\cd\w(\zeta-\del))$. Let
$I=\int_0^\cdot\bm \bm\cd dW$ and $I_k=\int_0^\cdot\bm_k\bm_k\cd
dW$. Then
$$
(I_k-I)(\cd\w(\zeta-\del))=\int_0^{\,\cd\,\w(\zeta-\del)}
 (\bm_k\bm_k^T-\bm\bm^T)dW,
$$
and using Burkholder's inequality and the convergence of $\bm_k$'s,
the left hand side of the last formula converges locally uniformly
to zero with probability one. Since $\del>0$ is arbitrary, we have
shown that the equation \eqref{11} holds for $Z$ and $\bm$ and all
$t<\zeta$. We would like $(X,Y)$ to satisfy the definition of a
mirror coupling given in Section \ref{secmirror} but at this point
we do not know whether $Y=X$ on $[\zeta,\infty)$. Hence, we redefine
$Y$ on $[\zeta,\infty)$ as $Y=X$ on this interval.
 The processes $X,Y,Z,L,M$ and $\bm$ and the random variable
$\zeta$ that we have constructed satisfy the definition of a
mirror coupling in $D$, given is Section \ref{secmirror}. As
follows from \cite{ab2}, the process $(X(t),Y(t))$ is strong
Markov. Lemma \ref{lemfour} applies to $(X,Y)$ because $(X,Y)$ is
a mirror coupling in $D$. Also \eqref{98} holds for $(X,Y)$ with
some $p_2>0$. It follows that the proof of Theorem \ref{th1}
presented above applies in the present setting, except for Lemma
\ref{lemf}, whose proof uses some properties of $\LL$. Hence, it
will suffice to prove that the assertion of Lemma \ref{lemf} holds
in the present context.

Part of the proof of Lemma \ref{lemf} uses some explicit
properties of $\LL$. It might be possible to derive the needed
properties of $\LL$ from those of $\LL_k$'s but that seems to be a
hard and tedious task so we will use an alternative approach.

Consider a second Neumann eigenfunction $\psi$ and assume that
$\psi\in S$. Consider any $(x,y)\in T$, $x\ne y$, and assume that
$(X(0),Y(0)) = (x,y)$. We have by \eqref{kb10},
$$
\psi(y)-\psi(x)=e^{\mu_2 t}\EE_{x,y}(\psi( Y(t))-\psi( X(t))).
$$
Since $\psi(x')\le\psi(y')$ for $(x',y')\in T$ and $(X(t),Y(t))\in
T$ for all $t<\zeta$, the right hand side is non-negative.
Moreover, the right hand side is strictly positive if for some
$t\geq 0$ we have $\PP_{x,y} (\psi( Y(t))>\psi( X(t)) >0$. Hence,
it remains to consider only the case when $\PP_{x,y} (\psi(
Y(t))>\psi( X(t)) =0$ for every $t$. By continuity of $X$ and $Y$,
this is equivalent to
 \be\label{kb30}
 \PP_{x,y} (\forall\ \psi( Y(t))=\psi( X(t)) =1.
 \ee
We have assumed that $x\ne y$ so $\zeta >0$, a.s. Reflected
 Brownian motion spends zero time on the boundary of $D$. Hence,
there exists a (random) time interval $[t_1,t_2]$ with $t_1<t_2<
\zeta$ such that $X(t),Y(t) \in D$ for all $t\in[t_1,t_2]$. It
follows that the mirror $\ell(t)$ does not move during this time
interval. Since $\psi$ is a real analytic function, we conclude
from this and \eqref{kb30} that $\psi$ is symmetric with respect
to $\ell^* := \ell(t_1)$.

Suppose that $D$ is symmetric with respect to $\ell^*$. Then
$\ell^*$ is either vertical or horizontal, by the assumption made
in the proposition. If $\ell^*$ is horizontal then $x$ and $y$ lie
on a vertical line but this is ruled out by the geometric
assumptions on $D_k$'s. If $\ell^*$ is vertical then $x$ and $y$
are symmetric with respect to the vertical line of symmetry of
$D$. Then $D_L$ and $D_R$ are also symmetric and it is easy to see
that the first part of the proof of Lemma \ref{lemf} applies and
one can conclude that $\psi(x)<\psi(y)$.

Next suppose that $D$ is not symmetric with respect to $\ell^*$.
Then there is a positive probability that one and only one of the
processes will spend some positive amount of local time on the
boundary of $D$. This will move the mirror before time $\zeta$ and
the same argument as before shows that there exists $\ell^{**}\ne
\ell^*$ that is a line of symmetry for $\psi$. Moreover,
$\ell^{**}$ can be chosen arbitrarily close to $\ell^*$. This
easily implies that either $\psi$ is constant, which is
impossible, or it is a function of only one variable in some
orthonormal coordinate system. An argument given in the proof of
Lemma 5 in \cite{ab2} shows that $D$ must be a rectangle. We have
assumed that $D$ is not a rectangle so the proof of the
proposition is complete.
 \qed

We believe that the assumptions on the domains $D_k$ converging to
$D$ eliminate the possibility that $D$ has a line of symmetry that
is not horizontal or vertical, or that $D$ is a rectangle, but
proving this does not seem to be useful. Hence, we added an
appropriate assumption about the shape of $D$ into Proposition
\ref{kb24}.

\section{Examples}\label{secex}

\beginsec

Most of the assumptions on $D$ listed in Section \ref{secassn}
must be checked directly in concrete examples; doing so is a
straightforward although tedious task. However, we would like to
comment on Assumption \ref{a4}. For any point $P\in \prt D$ with
$P_1 < P < P_3$, one can find $[P,Q]\in \calA(P_1,P_3)$ with
$\angle(P,Q)$ arbitrarily close to $\angle(P)$, by the definition
of $\calA(P_1,P_3)$. It is clear, therefore, that Assumption
\ref{a4} can be satisfied only if the curvature of $\prt D$ is
decreasing in a neighborhood of $P$. Vice versa, if the curvature
of $\prt D$ is non-increasing between $P_1$ and $P_3$ then the
assumption is satisfied for $[P,Q]\in \calA(P_1,P_3)$ with
$\angle(P,Q)$ very close to $\angle(P)$ or ``moderately'' close to
$\angle(P)$. For larger angles, Assumption \ref{a4} has to be
verified directly.

\begin{example}\label{ex1} {\rm
We will analyze the domain $D$ depicted in Fig.~\ref{figd}. The
following conditions uniquely define the domain.

\begin{itemize}

\item[1.] The domain $D$ is convex and its boundary passes through
points $(-2\sqrt{2},0)$, $(0,2)$, $(3,0)$ and $(0, -\sqrt{2})$.

\item[2.] The boundary is a piece of an ellipse between points
$(0,2)$ and $(3,0)$, with horizontal and vertical tangents at
endpoints.

\item[3.] The boundary is a piece of an ellipse between points
$(3,0)$ and $(0, -\sqrt{2})$, with horizontal and vertical
tangents at endpoints.

\item[4.] The boundary is a circular arc with center at $(0,-1)$
and endpoints $(-2\sqrt{2},0)$ and $(0,2)$. Note that the tangent
line is horizontal at $(0,2)$ but it is not vertical at
$(-2\sqrt{2},0)$.

\item[5.] The boundary is a piece of circular arc with center at
$(-\sqrt{2},0)$ and endpoints $(-2\sqrt{2},0)$ and $(-\sqrt{2}, -
\sqrt{2})$, with horizontal and vertical tangents at endpoints.

\item[6.] The boundary is a horizontal line segment between points
$(-\sqrt{2}, - \sqrt{2})$ and $(0, -\sqrt{2})$.

\end{itemize}

\begin{figure}
\centering
 \includegraphics[width=10cm]{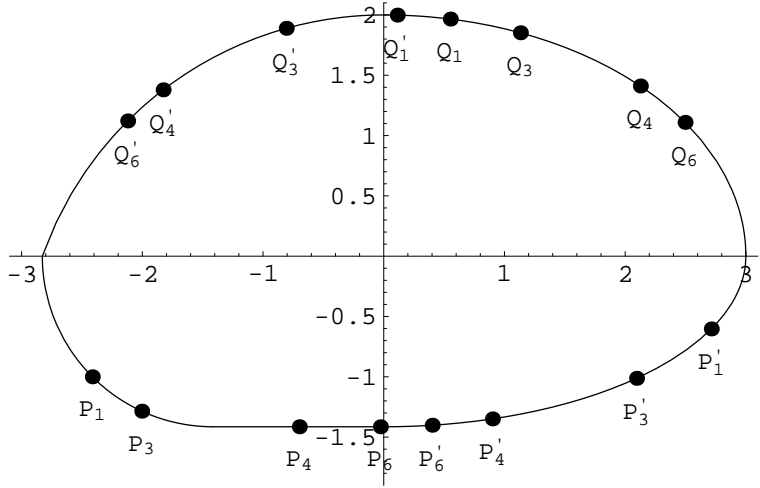}
\caption{A domain with simple second Neumann
eigenvalue.}\label{figd}
\end{figure}

The domain $D$ is not strictly convex and it is not $C^1$. We will
ignore these facts for the moment and we will proceed with a
choice of parameters and special points as in Section
\ref{secassn}. We take $\al=\pi/4$. This and the assumptions in
Section \ref{secassn} define uniquely points $P_1, P_3, P_4, P_6,
Q_1, Q_3, Q_4, Q_6$, and the analogous points with primes. We will
now describe how these points can be identified. $P_1$ is the
unique point with $\angle (P_1) = \pi/4$. $Q_1$ is the unique
point on the boundary with $\angle(P_1,Q_1) =\pi/4$. $Q'_6$ is the
point with $\angle(Q'_6)=7\pi/4$ and $P'_6$ is defined by
$\angle(P'_6,Q'_6) =7\pi/4$. The line segment $[P_3,Q_3]$ is
chosen so that $\angle(P_3,Q_3) =\pi/4$ and $[P_3,Q_3]\cap
[P'_6,Q'_6]$ is the midpoint of $[P'_6,Q'_6]$. Similarly,
$[P'_4,Q'_4]$ is chosen so that $\angle(P'_4,Q'_4) =7\pi/4$ and
$[P_1,Q_1]\cap [P'_4,Q'_4]$ is the midpoint of $[P_1,Q_1]$. Other
points are defined in the analogous way.

Because of the way the domain in our example is defined, the
coordinates of all special points are algebraic numbers and can be
written as explicit formulas involving only square roots.
 Some of the formulas are very complicated so we give coordinates
of the special points in the approximated decimal form. See also
Fig.~\ref{figd}.
\begin{align}\notag
P_1 &=(-2.41,-1.00), \ P_3=(-2.005,-1.28),\ P_4 = (-0.7,-1.41),
\ P_6=(-0.027,-1.41), \\
\notag Q_1 &=(0.55,1.97), \ Q_3= (1.13,1.85), Q_4= (2.13,1.41), \
Q_6=(2.50,1.11),\\
\notag P'_1 & = (2.71,-0.6),\ P'_3 = (2.09,-1.01), \ P'_4 =
(0.9,-1.35), \ P'_6 =(0.4,-1.40), \\
\notag Q'_1 &= (0.11,2), \ Q'_3 =(-0.81,1.89), \ Q'_4=
(-1.83,1.38), \ Q'_6= (-2.12,1.12).
\end{align}

We comment now on why $[P_3,Q_3]$ has been chosen so that
$[P_3,Q_3]\cap [P'_6,Q'_6]$ is the midpoint of $[P'_6,Q'_6]$. Note
that if $\ell(P',Q')$ is parallel to $[P_3,Q_3]$ and $P'> P_3$
then there are no lower left hinges for $\ell(P',Q')$. On the
other hand, if $P'=P$ and $Q'< Q_3$ then there exists a lower left
hinge for $\ell(P',Q')$. This observation is the basis of
verification of Assumptions \ref{a2} and \ref{a3}.

As for other assumptions listed in Section \ref{secassn}, some of
them are elementary but tedious to verify so we omit the formal
proof. The ones that are least trivial have been discussed at the
beginning of this section. Also, the assumptions of Proposition
\ref{kb24} regarding the domains $D_R$ and $D_L$ follow from
similar properties for the approximating sequence of domains.

Finally, note that because a part of $\pl D$ is a circular arc,
Assumption \ref{a1} does not hold for some line segments such that
 $\angle(P,Q)=\angle(P)$. We have to address this as well as the
fact that $D$ is not strictly convex and it is not $C^2$-smooth.
Approximating the circular arc by that of an ellipse, it is easy to
see that one can find a sequence of strictly convex $C^2$-smooth
domains $D_k \uparrow D$, where Assumption \ref{a1} holds for each
$D_k$. Moreover, $D_k$ can be chosen
 so that the points analogous to $P_j$'s, $Q_j$'s, $P'_j$'sand
$Q'_j$'s and defined relative to $D_k$ converge to the analogous
points in $D$. Hence, we can apply Proposition \ref{kb24} and we
conclude that the second Neumann eigenvalue in $D$ is simple.
Corollary \ref{kb23} implies that the second eigenfunction $\psi$
(or $-\psi$) is positive to the left of $[P_1,Q'_6]$ and negative
to the right of $[P'_1,Q_6]$. By Corollary \ref{kbgrad},
$\angle(\nabla \psi(x)) \in [0, \pi/4]\cup (3\pi/4, \pi)$ for
$x\in \ol D$ between $[P_3, Q'_4]$ and $[P'_3, Q_4]$.
 }
\end{example}

\begin{example}\label{ex2} {\rm
Our next example is related to \cite{jernad}, \cite{pas}
 and an earlier article \cite{banbur}.
Jerison and Nadirashvili proved in \cite{jernad} that the hot spots
conjecture holds in all convex planar domains with two perpendicular
axes of symmetry for all eigenfunctions corresponding to the second
eigenvalue but they left the question of the eigenvalue multiplicity
open. Pascu proved in \cite{pas} that the hot spots conjecture holds
for planar convex domains with a single line of symmetry, i.e., the
maximum and minimum of the second Neumann eigenfunction are attained
at the boundary of the domain. However, his theorem is stated for
only one of many possible eigenfunctions corresponding to the second
eigenvalue. The domain shown in Fig.~\ref{figf} has the boundary
consisting of two circular arcs and two line segments. Since the
ratio of its diameter to width is less than 1.53, Proposition 2.4 of
\cite{banbur} does not apply and we do not think that there is any
other theorem in the literature that implies that the second Neumann
eigenvalue is simple in this domain. This is indeed the case but we
omit the detailed proof as it follows the lines of Example
\ref{ex1}. We conclude that, in view of \cite{pas}, the hot spots
conjecture holds in its strongest form for the domain in
Fig.~\ref{figf} and similar convex symmetric planar domains with at
least one axis of symmetry.

\begin{figure}
\centering
 \includegraphics[width=6cm]{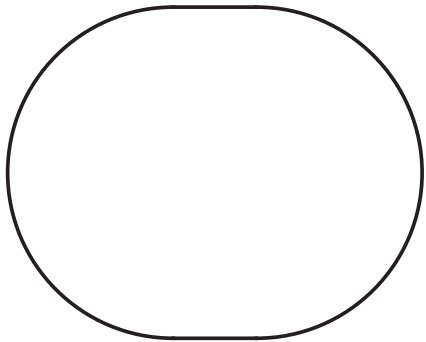}
\caption{A domain with a line of symmetry and small diameter to
width ratio.}\label{figf}
\end{figure}

 }
\end{example}

\begin{example}\label{ex3} {\rm
We conclude with a challenge for the reader, similar in spirit to
Exercise 4.1 in \cite{banbur}. That exercise is concerned with a
``snake'' domain, i.e., a twisted version of a very thin ``lip
domain,'' defined at the beginning of Section \ref{intro}. Our
present example is depicted in Fig.~\ref{figg}. One can show that
there exist subsets $D_L$ and $D_R$ of $D$ (close to the
``endpoints'' of $D$) and $T\subset \ol D \times \ol D$ such that
\eqref{99} holds. Then an argument similar to that in the proof of
Proposition \ref{kb24} can be used to show that the second Neumann
eigenvalue is simple in this domain. We leave the details of the
proof as an exercise for the reader.

\begin{figure}
\centering
 \includegraphics[width=6cm]{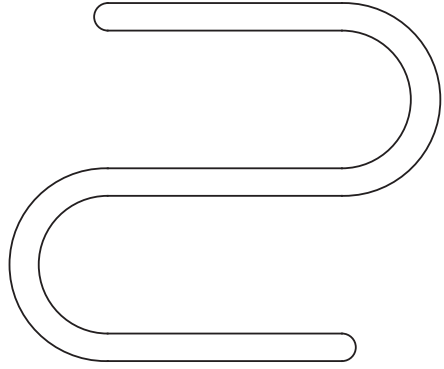}
\caption{A snake domain.}\label{figg}
\end{figure}

 }
\end{example}

\bibliographystyle{plain}

\end{document}